\documentclass[12pt,reqno]{amsart}
\usepackage{amssymb}
\usepackage{amsmath}

\usepackage{amssymb, amsmath, amsfonts, latexsym}
\usepackage{enumerate}

\newtheorem{thm}{Theorem}[section]
\newtheorem{cor}[thm]{Corollary}
\newtheorem{lem}[thm]{Lemma}
\newtheorem{prop}[thm]{Proposition}
\newtheorem{defn}[thm]{Definition}
\newtheorem{example}[thm]{Example}
\newtheorem{remarks}[thm]{Remark}

\numberwithin{equation}{section} \theoremstyle{remark}

\title[Ricci curvature and $W_1$-exponential convergence]{Ricci curvature and $W_1$-exponential convergence of Markov processes on graphs}

\author{Lingyan Cheng}
\address{Lingyan Cheng. Center for Applied Mathematics, Tianjin University, Tianjin 300072, PR China.}
\email{chengly@amss.ac.cn}

\author{Ruinan Li}
\address{Ruinan Li. School of Statistics and Information, Shanghai University of International Business and Economics, Shanghai 201620,  PR China.}
\email{ruinanli@amss.ac.cn}

\author{Liming Wu}
\address{Liming Wu. Laboratoire de Math\'ematiques Blaise Pascal, CNRS-UMR
6620, Universit\'e Clermont-Auvergne (UCA), Campus Universitaire des Cezeaux, 3 Place Vasarely, 63178 Aubi\`ere, France.}
\email{Li-Ming.Wu@math.univ-bpclermont.fr}

\newcommand{\cc}{\mathbb{C}}

\newcommand{\ee}{\mathbb{E}}

\newcommand{\nn}{\mathbb{N}}
\newcommand{\rr}{\mathbb{R}}
\newcommand{\pp}{\mathbb{P}}

\newcommand{\zz}{\mathbb{Z}}

\def\BB{\mathcal B}

\def\FF{\mathcal F}
\def\GG{\mathcal G}

\def\LL{\mathcal L}
\def\MM{\mathcal M}

\def\vep{\varepsilon}
\def\<{\langle}
\def\>{\rangle}

\def\sgn{{\rm sgn}}

\def\beq{\begin{equation}}
\def\nneq{\end{equation}}

\def\bdef{\begin{defn}}
\def\ndef{\end{defn}}

\def\bthm{\begin{thm}}
\def\nthm{\end{thm}}

\def\bprop{\begin{prop}}
\def\nprop{\end{prop}}

\def\brmk{\begin{remarks}}
\def\nrmk{\end{remarks}}

\def\bexa{\begin{example}}
\def\nexa{\end{example}}

\def\blem{\begin{lem}}
\def\nlem{\end{lem}}

\def\bcor{\begin{cor}}
\def\ncor{\end{cor}}

\date{}

\def\bexe{\begin{exe}}
\def\nexe{\end{exe}}

\def\bprf{\begin{proof}}
\def\nprf{\end{proof}}

\def\benu{\begin{enumerate}}
\def\nenu{\end{enumerate}}

\def\bdes{\begin{description}}
\def\ndes{\end{description}}

\begin{document}
\maketitle

\begin{abstract} In this paper, we show that the Ricci curvature lower
bound in Ollivier's Wasserstein metric sense of a continuous time jumping
Markov process on a graph can be characterized by some optimal coupling  generator and provide the construction of this latter. Some previous results of Ollivier for discrete time Markov chains are generalized to the actual continuous time case.
We propose a comparison technique with some death-birth
process on $\nn$ to obtain some explicit exponential convergence rate, by modifying the metric. A counterpart of Zhong-Yang's estimate is established in the case where the Ricci curvature with repsect to the graph metric is nonnegative.
Moreover we show that the Lyapunov function method for the exponential convergence works with some explicit quantitative estimates, once if the Ricci curvature is bounded from below by a negative constant. Finally we present applications to Glauder dynamics under some dynamical versions of the Dobrushin uniqueness condition or of the Dobrushin-Shlosman analyticity condition.

\end{abstract}



\section {Introduction}

\subsection{Graph, generator, Ricci curvature lower bound}
Let $G=(S, E)$ be an at most countable connected graph with vertex
set $S$ and oriented edges set $E$, which is a symmetric subset of
$S^2\backslash \{(x,x);x\in S\}$. If $(x,y)\in E$, we call that
$x,y$ are neighbors or adjacent, denoted by $x\sim y$. We assume
always that the degree $d_x=|\{y\in S: y\sim x\}|$ (for a set $A$,
$|A|$ denotes the number of the elements in $A$) is finite for every
$x\in S$, that is, $G$ is locally finite. A path $\gamma_{xy}$ from
$x$ to $y$ is a family of edges $\{e_1,\cdots, e_n\}$ where
$e_k=(x_{k-1},x_k)\in E$, such that $x_0=x,x_n=y$, and its graphic
length is defined as the number $n$ of edges in this path, denoted
by $|\gamma_{xy}|_G$. A geodesic from $x$ to $y$ ($x\ne y$) is a
path from $x$ to $y$ with a minimal number of edges, and its graphic
length is defined as the graph distance $d_G(x,y)$ between $x$ and
$y$. Of course $d_G(x,x):=0$. A great difference of graphs from
Riemannian manifolds is that the geodesics between two vertices on
graphs are in general not unique except trees.

Consider a Markov  generator of nearest-neighbor type:
\beq\label{LL} \LL f(x)=\sum_{y \in S} J(x,y)\left(f(y)-f(x)\right), \ \
\text{for all} \ x \in S \nneq for any function $f: S\rightarrow
\rr$, where the jump rate $J(x,y)$ from $x$ to $y$ is non-negative,
and $J(x,y)>0$ if $x\sim y$ and $J(x,y)=0$ if $d_G(x,y)\ge 2$. When
$J(x,y)=\frac 1{d_x}$, $\LL$ is the {\it Laplacian operator}
$\Delta$ on $G$.

Let $(X_t)$ be the Markov process generated by $\LL$, defined on
$(\Omega, (\FF_t)_{t\ge0}, (\pp_x)_{x\in S})$, with the transition
probability semigroup $P_t=e^{t\LL}$ generated by $\LL$. We assume
that it is conservative and it has an invariant probability measure
$\mu$ (necessarily unique). It can be viewed as a Monte-Carlo
algorithm for sampling $\mu$. One basic question is the exponential
convergence rate of $P_t$ to $\mu$ at large time. The reader is
referred to the books:  M.F. Chen \cite{Chen92, Chen05} (for the
methods of coupling and of functional inequalities), Meyn and
Treedie \cite{MT} (method of Lyapunov functions), and L.
Saloff-Coste \cite{SaC97} (functional inequalities).

A powerful approach to this question of exponential convergence is
to use optimal transport (or coupling in probability language).
Given a cost function $c(x,y)$ on $S^2$ such that $c(x,y)=c(y,x)>0$
for $x\ne y$ and $c(x,x)=0$. The transport cost $T_c(\nu,\mu)$ from
a probability measure $\nu$ to $\mu$ is defined as \beq\label{TC}
T_c(\nu,\mu)=\inf_{P\in \Pi(\nu,\mu)}\iint_{S^2} c(x,y) dP(x,y),
\nneq where $\Pi(\nu,\mu)$ are the family of all couplings of
$(\nu,\mu)$, i.e. probability measures on $S^2$ with marginal
distributions $\nu,\mu$. This is the probabilistic reformulation of
the Monge-Ampere optimal transport problem by Kantorovitch.

When $c(x,y)$ is some metric $d(x,y)$ on $S$, the transport cost
$T_c(\nu,\mu)$ becomes the $L^1$-Wasserstein distance
$W_{1,d}(\nu,\mu)$. We are interested in the estimate of constants
$\kappa>0$ and $K\ge 1$ such that \beq\label{expC} W_{1,d}
(P_t(x,\cdot), P_t(y,\cdot))\le K e^{-\kappa t} d(x,y), \ \forall
t\ge0, x,y\in S.\nneq If $S$ is finite, let $\lambda_1$ be the first
eigenvalue of $-\LL$ with the smallest real part $Re(\lambda_1)$,
then the best possible constant $\kappa$ is $Re(\lambda_1)$, but the
estimate of $K$ with $\kappa=Re(\lambda_1)$ becomes very difficult
or even impossible if the algebraic multiplicity of $\lambda_1$ is
strictly greater than the geometric one. A general idea to obtain
global behavior (for all time $t$) in (\ref{expC}) is to find a new
metric  $\tilde d$,  equivalent to $d$ (i.e. $M^{-1}d\le \tilde d\le
M d$ for some constant $M\ge1$), so that (\ref{expC}) holds with
$K=1$ w.r.t. $\tilde d$.

When $K=1$, (\ref{expC}) is equivalent to the following infinitesimal time
version \beq\label{Ric} \limsup_{t\to0+}\frac{ W_{1,d}
(P_t(x,\cdot), P_t(y,\cdot))-d(x,y)}{t}\le -\kappa d(x,y), \ \forall
x,y\in S. \nneq If it holds we say that {\it the Ricci curvature of
$\LL$ w.r.t. the metric $d$ is bounded from below by $\kappa$
($\in\rr$)}, written as $Ric(\LL, d)\ge \kappa$. When  $J(x,y)$ is a
transition probability kernel (i.e. $\sum_y J(x,y)=1$ for every
$x\in S$), this definition of Ricci curvature lower bound goes back
to Ollivier \cite{Oll09} in which the Ricci curvature of  discrete
time Markov chains is introduced and studied. When $\LL=\Delta$, the Laplacian on a connected Riemannian manifold equipped with the Riemannian
metric $d$, Renesse and Sturm \cite{RSturm} proved that $\kappa$ coincides with the
lower bound of the Ricci curvature $Ric$, showing that
the above definition is a natural one.

If (\ref{Ric}) holds for the Laplacian $\LL=\Delta$ w.r.t. the graph
metric $d=d_G$, we say that {\it the Ricci curvature of the graph
$G$ is bounded from below by $\kappa$, written as $Ric^G\ge \kappa$
} simply. See Lin, Lu and Yau \cite{LLY11} for some earlier
definitions of the Ricci curvature lower bound and studies.

As $P_t$ is unknown, a usual way in probability to obtain
(\ref{Ric}) is to construct a coupling Markov generator $\LL^\pi$ on
$S^2$ of $\LL$, i.e. $\LL^\pi (f\oplus g)(x,y)=\LL f(x) + \LL g(y)$
for all $(x,y)\in S^2$ where $(f\oplus g)(x,y)=f(x)+g(y)$, such that
\beq\label{RicL} (\LL^\pi d)(x,y)\le -\kappa d(x,y). \nneq

\subsection{$J$ is a transition probability matrix}
When $J(x,y)$ is a transition probability kernel (i.e. $\sum_y
J(x,y)=1$ for every $x\in S$), the process generated by $\LL$ is
given by $X_t=Z_{N(t)}$ where $(Z_n)_{n\in\nn}$ is the Markov chain
with transition probability kernel $J$ and $N(t)$ is a Poisson
process with parameter $\lambda=1$, independent of $(Z_n)$.
Following Ollivier \cite{Oll09}, {\it the Ricci curvature of $J$ is
said to be bounded from below by some constant $\kappa$, if
\beq\label{Ollivier1} W_{1,d}(J(x,\cdot), J(y,\cdot)) - d(x,y)\le
-\kappa d(x,y), \ x,y\in S. \nneq} For every $(x,y)\in S^2$, let
$J^\pi((x,y),\cdot)$ be an optimal coupling of $J(x,\cdot),
J(y,\cdot)$ realizing the Wasserstein metric $W_{1,d}(J(x,\cdot),
J(y,\cdot))$, then the corresponding generator with jumping  rates
$J^\pi$
$$
\LL^\pi F(x,y) =\sum_{(x',y')} (F(x',y')-F(x,y)) J^{\pi}((x,y),
(x',y'))
$$
satisfies $\LL^{\pi} d (x,y)\le -\kappa d(x,y)$. Thus $Ric(\LL,
d)\ge \kappa$.

Take a simple example to show that this approach from the definition
of Ricci curvature of the discrete time transition probability
kernel $J$ to that of the continuous time generator $\LL=J-I$ is far
from being sharp: let $S=\{0,1\}$ and $J(0,1)=J(1,0)=1$, $\LL$ is
the Laplacian. It is easy to see that $Ric(\Delta, d_G)=2$, however
the Ricci curvature of $J$  is $0$.

That is why Lin, Lu and Yau \cite{LLY11} introduced for every
$\alpha\ge 0$, {\it $\alpha$-Ricci curvature lower bound
$\kappa^{(\alpha)}$ } as
$$
(1+\alpha) [W_{1,d}(\tilde J_x, \tilde J_y)-d(x,y)]\le
-\kappa^{(\alpha)} d(x,y)
$$
where $\tilde J_x(z)= \frac 1{1+\alpha} (J(x,z) + \alpha
\delta_x(z))$ ($\delta_x$ is the Dirac measure at $x$). If $\tilde
J^\pi((x,y), \cdot)$ is an optimal coupling of $\tilde J_x$ and
$\tilde J_y$, then the corresponding generator $\LL^\pi$ with
jumping rates $(1+\alpha)\tilde J^\pi$ is a coupling of $\LL$ and
satisfies $\LL^\pi d\le - \kappa^{(\alpha)} d$. Thus $Ric(\LL, d)\ge
\kappa^{(\alpha)}$.

\subsection{Objective of the paper.}

Our first natural question is whether the Ricci curvature lower
bound $\kappa$ defined in (\ref{Ric}) is achieved by a coupling
$\LL^{\pi}$ of $\LL$ satisfying (\ref{RicL}). We will show that this
is true, providing so a generator's criterion for the Ricci curvature
lower bound.

The next question is how to construct an optimal coupling $\LL^\pi$
of $\LL$ in the sense that $\LL^{\pi}d(x,y)=\frac {d}{dt}
W_{1,d}(P_t(x,\cdot), P_t(y,\cdot))|_{t=0}$ for every $(x,y)\in
S^2$.

Our third question is: when $Ric(\LL, d_G)$ is bounded from below
by a non-positive constant, could one find some other metric $d$ so
that $Ric(\LL, d)\ge \kappa>0$ ? if yes, how ?

We will answer the first two questions in a positive way in full generality for the Markov processes on graphs. Answer to the third question above depends
certainly on studied models: we study it through comparison with a
birth-death process on $\nn$ and Lyapunov test-function method.

\subsection{Organization of the paper}
In the next section, we will show that the Ricci curvature can be characterized by some optimal coupling generator $\LL^\pi$,  furnishing so an infinitesimal criterion for the Ricci curvature. Moreover we provide the construction of the optimal coupling generator(s). In Section 3, we introduce some comparison technique and the comparison condition $\bf C(J, \alpha,\beta)$ with one dimensional birth-death processes, which allow us to improve Lin-Yau's estimate of the spectral gap in terms of the degree and the diameter. We establish  the exponential convergence of $(X_t)$ in $W_{1,d_G}$ even if the Ricci curvature w.r.t. the graph metric $d_G$ is bounded below by a negative constant, when it is dissipative at infinity (a generalization of Eberle's \cite{Eberle16} result for diffusions to graphs). A counterpart of the famous Zhong-Yang's estimate for the spectral gap $\lambda_1$ when the Ricci curvature w.r.t. the graph metric $d_G$ is non-negative is established in Section 5. In Section 6, we propose a very practical criterion of Lyapunov function type for the exponential convergence in $W_1$, generalizing the result of Hairer-Mattingly \cite{HM11} from the discrete time Markov chains to continuous time jumps processes. As application the last section is devoted to the exponential convergence in $W_1$ for the Glauber dynamics on high dimensional product graphs, under the Dobrushin's uniqueness condition or Dobrushin-Shlosman analyticity condition.

\section{Optimal coupling of generator}
Given the generator $\LL$ in (\ref{LL}), the Markov process $(X_t)$
can be easily realized by stochastic algorithm: if $X_0=x$, then $X_t=x$ until
the first jumping  time $\tau_1=\xi_1/\lambda(x)$ where $\xi_1$ is a
random variable of exponential law with parameter $1$, and
\beq\label{jumps-rate} \lambda(x)=\sum_{x'\ne x} J(x,x') \nneq is
the total jumping  rate at $x$; and at the first jumping time $\tau_1$, $X_{\tau_1}$ is distributed as
$J(x,\cdot)/\lambda(x)$. As the (conditional) distribution of
$X_{\tau_1}$ does not depend upon $\tau_1$, $\tau_1$ and
$X_{\tau_1}$ are independent. Next run the algorithm by regarding
$(X_{\tau_1}, \tau_1)$ as a new starting space-time point,
and so on for obtaining the {\it n-th} jumping  time $\tau_n$ and
$X_{\tau_n}$ $(n\in \nn)$.

But for the coupling purpose we may add an artificial {\it false}
jumping  rate $J(x,x)\ge0$, which does not change the generator
$\LL$. In the algorithm above, $\tau_1$ becomes now $\xi_1/J(x,S)$
and $X_{\tau_1}$ is distributed as $J(x,\cdot)/J(x,S)$. The only
difference is that the probability of the false jump
$\pp_x(X_{\tau_1}=x)=\frac{J(x,x)}{J(x,S)}$ may be positive.

Since $M_t(f)=f(X_t)-f(X_0)-\int_0^t \LL f(X_s)ds$ is a local
martingale and $\sup_{t\le \tau_1}|M_t(f)|\le |f(x)| +
|f(X_{\tau_1})|+\tau_1 |\LL f(x)|$, we have $\ee_x M_{\tau_1}(f)=0$
and then
$$
\LL f(x) = \frac {\ee_x f(X_{\tau_1}) - f(x) }{\ee \tau_1}.
$$
This can be verified directly.

As said in the introduction we assume always that the process is
conservative and has an invariant probability measure $\mu$
(necessarily unique as our graph is connected and $J(x,y)>0$ if
$y\sim x$). The conservativeness and the existence of $\mu$  are
equivalent to the positive recurrence, which is characterized by the
following well known Lyapunov function criterion (see Meyn and
Tweedie \cite{MT}):

{\bf (H)} {\it There are two positive functions $V\ge 1$ and $U$ on
$S$ with $\inf_S U>0$,  and some finite subset $K$ of $S$, and a
positive constant $b$ such that \beq\label{Lyapunov1} \LL V \le -U +
b1_K. \nneq} Indeed if $U=1$, that is the positive recurrence; if
$U=-\delta V$  for some positive contant $\delta$, that is
equivalent to the exponential recurrence. This hypothesis is and will be assumed throughout the paper.

\subsection{Coupling of generators}

\bdef {\it Given two Markov generators $\LL_1$ and $\LL_2$ on $S$,
with jumping  rates kernels $J_1$ and $J_2$, a coupling Markov
generator $\LL^\pi$ of $\LL_1$ and $\LL_2$
 is defined as
$$
\LL^\pi F(x,y) := \sum_{(x',y')\in S^2} [F(x',y')-F(x,y)]
J^\pi((x,y), (x',y')),\ (x,y)\in S^2
$$
for all functions $F:S^2\to \rr$ with finite support, where $J^\pi$
is a nonnegative kernel on $S^2$, such that for any functions $f,g$
on $S$ with finite support,
$$
\LL^{\pi}(f\oplus g)(x,y) = \LL_1 f(x) + \LL_2 g(y), \ (x,y)\in S^2
$$
where $(f\oplus g)(x,y):=f(x) +g(y), (x,y)\in S^2$.
That is equivalent to: for every $(x,y)\in S^2$,
 \beq\label{21}\aligned &\sum_{y'}
J^\pi((x,y), (x',y')) = J_1(x,x'), \ \forall x'\ne x;\\
& \sum_{x'} J^\pi((x,y), (x',y')) = J_2(y,y'),\ \forall y'\ne y.
\endaligned\nneq
 A coupling Markov generator $\LL^\pi$ of $\LL$ is that of
$\LL_1=\LL_2=\LL$, satisfying moreover \beq\label{22} J^\pi((x,x),
(x',x'))=J(x,x'), \ J^\pi((x,x), S^2\backslash \Delta)=0, \nneq
where
$\triangle=\{(x,x); \ x\in S\}$ is the diagonal of $S^2$. The
jumping rates kernel $J^\pi$ of $\LL^\pi$ will be called coupling
jumping rates kernel of $J$.} \ndef

The last condition (\ref{22}) means that the Markov process
$(X_t,Y_t)$ generated by the coupling generator $\LL^\pi$ of $\LL$,
once getting together, will be together forever.  The process $(X_t,Y_t)$
generated by any coupling generator $\LL^\pi$ of $\LL$
is also conservative, since for any finite subset $K$ of $S$, the
exiting time $\sigma_{K^2}(X,Y)=\inf\{t\ge 0; (X_t,Y_t)\notin K^2\}$
of $(X_t,Y_t)$ is equal to the minimum $\sigma_K(X)\wedge
\sigma_K(Y)$ of the exiting times of $X$ and $Y$ from $K$, and
$\sigma_K(X)\wedge \sigma_K(Y)\to +\infty$ as $K\uparrow S$ by the
conservativeness of $X$ and $Y$, assured by our assumption {\bf
(H)}.

A natural coupling generator of $\LL_1$ and $\LL_2$ in probability
is {\it the independent coupling} whose joint jumping rates kernel
is given by
\beq\label{indepC} J^\pi((x,y), (x',y'))
=1_{y'=y}J_1(x, x') + 1_{x'=x} J_2(y, y'),
\nneq
for
all $(x',y')\ne (x,y)$ in $S^2$. When this holds at a fixed couple
$(x,y)\in S^2\backslash \triangle$ for a coupling generator
$\LL^\pi$ of a single generator $\LL$, as
$$\aligned
\frac d{dt}&[P_t^\pi (f\otimes g)(x,y) - P_tf(x)
P_tg(y)]|_{t=0}\\
&=\LL^\pi (f\otimes g)(x,y) - [g(y)\LL f(x)
+f(x)\LL g(y)]\\
&=\sum_{(x',y')\in S^2} (f(x')-f(x))(g(y')-g(y))J^\pi((x,y),
(x',y'))=0
\endaligned$$ for any bounded functions $f,g$ on $S$, where
$(f\otimes g)(x,y):=f(x)g(y)$, $X_t$ and $Y_t$ are asymptotically
independent for small $t$. In that case we say that $\LL^\pi$ is
{\it a locally independent coupling} of $\LL$ at $(x,y)$.  If
(\ref{indepC}) holds for $J_1=J_2=J$ and for all $(x,y)\in
S^2\backslash \triangle$, we say again, with some abuse, $\LL^\pi$
is {\it the independent coupling} of $\LL$, though $X_t=Y_t$ after
the coupling time $\tau_c:=\inf\{s\ge0; X_s=Y_s\}$ and they are not
independent after $\tau_c$.

 The following theorem says that the Ricci curture
lower bound can be characterized by coupling generator, answering
the first question raised in the Introduction.

\bthm\label{thm21} Let $c(x,y)$ be a cost function (i.e.
$c(x,y)=c(y,x)>0$ and $c(x,x)=0$ for all $(x,y)\in S^2\backslash
\triangle$).

\benu[(a)] \item There is always a $T_c$-optimal coupling generator
$\LL^\pi$ of $\LL$ in the sense that \beq\label{thm21c}
\liminf_{t\to 0} \frac{T_c(P_t(x,\cdot),P_t(y,\cdot))-c(x,y)}{t} \ge
\LL^\pi c(x,y), \ \forall (x,y)\in S^2. \nneq

\item
Given some constant $\kappa\in \rr$, \beq\label{thm21a}
T_c(P_t(x,\cdot),P_t(y,\cdot))\le e^{-\kappa t} c(x,y), \ \forall
t\ge0, (x,y)\in S^2 \nneq holds if and only if there is some
coupling generator $\LL^\pi$ of $\LL$ such that \beq\label{thm21b}
\LL^\pi c(x,y) \le -\kappa c(x,y), \ \forall (x,y)\in S^2\backslash
\triangle. \nneq \nenu In particular if $c(x,y)=d(x,y)$ is some
metric on $S$, $Ric(\LL,d)\ge \kappa$ if and only if (\ref{thm21b})
holds with $c(x,y)=d(x,y)$ for some coupling generator $\LL^\pi$.

 \nthm

\bprf {\bf (a).} For every $h>0$ and $(x,y)\in S^2$, let
$P_{(h)}^\pi((x,y), \cdot)$ be an optimal coupling of $P_h(x,\cdot),
P_h(y,\cdot)$ such that
$$
\sum_{(x',y')\in S^2} c(x',y') P_{(h)}^\pi((x,y), (x',y')) =
T_c(P_h(x,\cdot),P_h(y,\cdot)).
$$
Necessarily  $P_{(h)}^\pi((x,x), S^2\backslash \Delta)=0$ and
$P_{(h)}^\pi((x,x), (x',x'))=P_h(x,x')$. Notice that
$(P_{(h)}^\pi((x,y), \cdot))_{h\ge0}$ is in general not a semigroup.

Let $\LL^\pi_{(h)}F=\frac 1h (P^\pi_{(h)}F -F)$ ($F$ is of finite
support) which is a Markov generator with jumping  rate kernel
$J^\pi_{(h)}=(1/h)P^\pi_{(h)}$ . If $F(x,y)=f(x) +g(y)$ where $f,g$
are functions on $S$ of finite support,
$$
\LL_{(h)}^\pi F(x,y)=\LL_{h}f(x) + \LL_h g(y)
$$
where $\LL_hf(x)=\frac 1h[P_hf(x)-f(x)]$. Remark that
$$
\lim_{h\to 0} \LL_{h} f(x)=\LL f(x)
$$
for any bounded function $f$ on $S$ such that $\LL f$ is bounded.

For every $(x_0,y_0)\in S$, letting $f_{x_0}(x)=1_{x\ne x_0}$, we
have
$$\aligned
J^\pi_{(h)}((x_0,y_0), S^2\backslash \{(x_0,y_0)\})&\le
J^\pi_{(h)}(f_{x_0}\oplus f_{y_0})(x_0,y_0)\\
&= \LL_{(h)}^\pi (f_{x_0}\oplus f_{y_0})(x_0,y_0) \\
&= \LL_h f_{x_0}(x_0) + \LL_h f_{y_0}(y_0) \\
&\to \LL f_{x_0}(x_0) + \LL f_{y_0}(y_0)=\lambda(x_0)+\lambda(y_0)
\endaligned$$
when $h\to0+$. Thus we can take a decreasing sequence $h_n\to 0$
such that \benu
\item for all different points $(x,y), (x',y')\in S^2$,
$$ J^\pi((x,y), (x',y')):=\lim_{n\to \infty}J^\pi_{(h_n)}((x,y), (x',y'))$$
exists in $\rr^+$;

\item for all $(x,y)\in S^2$,
$$\aligned
&\liminf_{t\to 0} \frac{T_c(P_t(x,\cdot),P_t(y,\cdot))-c(x,y)}{t}\\
& = \lim_{n\to \infty} \frac 1{h_n}(T_c(P_{h_n}(x,\cdot),
P_{h_n}(y,\cdot))-c(x,y))\in \overline{\rr}:=[-\infty, +\infty].
\endaligned$$
 \nenu
Define
$$\LL^\pi F(x,y) = \sum_{(x',y')\ne (x,y)}
[F(x',y')-F(x,y)]J^{\pi}((x,y), (x',y')), \ (x,y)\in S^2.$$ It is a
coupling Markov generator of $\LL$.

Now let $(S_m)_{m\ge1}$ be an increasing sequence of finite subsets such
that $\bigcup_{m \ge 1} S_m=S$. For any $(x,y)\in S^2$ and for $m$ large
enough, $(x,y)\in S_m^2$, and
$$\aligned
&\frac 1{h_n}[T_c(P_{h_n}(x,\cdot), P_{h_n}(y,\cdot))-c(x,y)] \\
&\ge \frac 1{h_n} \sum_{(x',y')\in S_m^2} P^\pi_{(h_n)}((x,y),
(x',y'))[c(x',y')-c(x,y)]\\
& \to \sum_{(x',y')\in S_m^2} J^\pi((x,y),
(x',y'))[c(x',y')-c(x,y)]\endaligned
$$
as $n\to +\infty$. That implies (\ref{thm21c}) by letting
$m\to+\infty$.

\medskip
{\bf (b).} The necessity follows by part (a). For the sufficiency,
let $(X_t,Y_t)$ be the Markov process generated by $\LL^\pi$. By
(\ref{thm21b}) and Ito's formula and Fatou's lemma, $e^{\kappa t}
c(X_t,Y_t)$ is a $\pp_{(x,y)}$-supermatingale. Thus
$$
T_c(P_t(x,\cdot),P_t(y,\cdot))\le \ee_{(x,y)}c(X_t,Y_t)\le
e^{-\kappa t} c(x,y)
$$
as desired.
 \nprf

\brmk{\rm Our argument in Theorem \ref{thm21}(a) (which implies the necessity of Theorem \ref{thm21}(b))
depends on the structure of jump processes valued in countable states space. How to extend it to diffusions is a quite interesting open question.
}\nrmk

\subsection{Construction of an optimal coupling generator}
The construction of an optimal coupling $\LL^\pi$ above, based on
the coupling of the unknown semigroup $P_t$, is only theoretical. We
turn now to the construction of $\LL^\pi$ directly from the jumping
rates kernel $J$. To this purpose we will at first extend the notion
of coupling of probability measures to any two nonnegative measures
of same mass.

Given two positive measures $\nu_1, \nu_2$ on $S$, a coupling of
$\nu_1, \nu_2$ is a positive measure $\nu^\pi$ on $S^2$ such that
$$
\sum_{y} \nu^\pi(x,y)=\nu_1(x), \ \sum_{x} \nu^\pi(x,y)=\nu_2(y), \
\forall x,y\in S
$$
i.e. its marginal measures are $\nu_1$ and $\nu_2$. Necessarily
$\nu_1(S)=\nu_2(S)=\nu^\pi(S^2)$: they must have the same mass. One
can define the transport cost $T_c(\nu_1,\nu_2)$ by (\ref{TC}). The
optimal couplings exist always for $T_c(\nu_1,\nu_2)$, once if it is
finite.

Now if $\LL^\pi$ is a coupling generator, by (\ref{21}) its jumping
rates measure $J^\pi((x,y),\cdot)$ will be a coupling of
$J(x,\cdot)$ and $J(y,\cdot)$ but with some false jump rates
$$
J(x,x)=\sum_{y'} J^\pi((x,y), (x,y')),\ J(y,y)=\sum_{x'} J^\pi((x,y),
(x',y)).
$$
Throughout this paper $J(x,\cdot)$ will be a nonnnegative measure on
$S$ with $J(x,x')$ being the jumping  rate of $\LL$ for $x'\ne x$,
and with some specified value of $J(x,x)$ (false jump  rate).

\bthm\label{thm22} For any $x\ne y$ and {\bf every} $\lambda(x,y)\ge
\lambda(x)+\lambda(y)$, let
$$
J^\pi((x,y),\cdot)
$$
be an optimal coupling of the following two  measures of mass
$\lambda(x,y)$: \beq\label{thm22a}\aligned
J(x, x') &= 1_{x'\ne x}J(x,x')+[\lambda(x,y)-\lambda(x)] \delta_x(x'),\\
J(y, y') &=1_{y'\ne y}
J(y,y')+[\lambda(x,y)-\lambda(y)]\delta_y(y')
\endaligned\nneq
realizing $T_c(J(x, \cdot), J(y,\cdot))$, where $\delta_x$ is the Dirac measure at point $x$. Then the generator
$\LL^\pi$ with jumping  rates kernel $ J^\pi((x,y), \cdot)$ is an
optimal coupling generator of $\LL$ in the sense that for every
$(x,y)\in S^2$ fixed, $\LL^\pi c(x,y)$ attains the minimum among all
coupling Markov generators.
\nthm

In practice one takes $\lambda(x,y)=\lambda(x)+\lambda(y)$ for all
$x\ne y$.

\bprf For any coupling Markov generator $\LL^\pi$ with jumping rates
kernel $J^\pi$, its total jumping  rate
$$\aligned
\lambda^\pi(x,y)&=\sum_{(x',y')\ne (x,y)} J^\pi((x,y),(x',y'))\\
&= \sum_{x'\ne x}J(x, x') + \sum_{y'\ne y}J^\pi((x,y),(x,y'))\in
[\lambda(x), \lambda(x)+\lambda(y)] \endaligned$$ by (\ref{21}).
Then $\max\{\lambda(x),\lambda(y)\}\le \lambda^\pi(x,y)\le
\lambda(x)+\lambda(y)$. For $(x_0,y_0)$ fixed, as the set of the
marginal measures of $x'$  of $J^\pi((x_0,y_0),\cdot)$ with $J^\pi$
varying over all coupling jumping  rates such that $J^\pi((x_0,y_0),
(x_0,y_0))=0$ is contained in $\{J(x_0,\cdot) + \beta \delta_{x_0};
\beta\in [0,\lambda(y_0)]\}$, it is tight. The same for the set of
marginal measures of $y'$ of $J^\pi((x_0,y_0),\cdot)$. Thus the set of all coupling rates
measures $J^\pi((x_0,y_0),\cdot)$ (i.e. satisfying (\ref{21})) with
$J^\pi((x_0,y_0), (x_0,y_0))=0$, is tight.

Therefore
$$
\inf_{\LL^\pi} \LL^\pi c(x_0,y_0) = \inf_{J^\pi} \sum_{(x',y')\in
S^2} [c(x',y')-c(x_0,y_0)] J^\pi((x_0,y_0), (x',y'))
$$
is attained at some jumping  rate measure $J^{\pi_0}((x_0,y_0),
\cdot)$, for the last functional is lower semi-continuous on $J^\pi$
w.r.t. the weak convergence topology. But
$$\aligned
&\LL^{\pi_0}c(x_0,y_0)=\sum_{(x',y')\in S^2} [c(x',y')-c(x_0,y_0)]
J^{\pi_0}((x_0,y_0),(x',y')).\\
\endaligned$$
By the expression above, with the
total jumping  rate $\lambda^{\pi_0}(x_0,y_0)$ fixed,  $
J^{\pi_0}((x_0,y_0),\cdot)$ must be an optimal coupling of its two
marginal measures
$$\aligned
J^{\pi_0}_x(x_0, x')= 1_{x'\ne x_0}J(x_0,x') +
[\lambda^{\pi_0}(x_0,y_0)-\lambda(x_0)]\delta_{x_0}(x'),\\
J^{\pi_0}_y(y_0, y')= 1_{y'\ne y_0}J(y_0,y') +
[\lambda^{\pi_0}(x_0,y_0)-\lambda(y_0)]\delta_{y_0}(y')\\
\endaligned$$ for the transport cost $T_c(J^{\pi_0}_x(x_0,
\cdot), J^{\pi_0}_y(y_0, \cdot))$. As $(x_0,y_0)$ is arbitrary, we
get an optimal coupling jumping  rates kernel $J^{\pi_0}$.

Now we show that for any jumping  rate $\lambda(x_0,y_0)\ge
\lambda^{\pi_0}(x_0,y_0)$, we can construct another optimal
coupling. Let $J^{\pi,\lambda}((x_0,y_0),\cdot)$ be an optimal
coupling of
$$\aligned
J^{\lambda}_x(x_0, x') &=1_{x'\ne x_0} J(x_0,x') +
[\lambda(x_0,y_0)-\lambda(x_0)]\delta_{x_0}(x')\\
&= J^{\pi_0}_x(x_0,x') +
[\lambda(x_0,y_0)-\lambda^{\pi_0}(x_0,y_0)] \delta_{x_0}(x'),\\
J^{\lambda}_y(y_0, y') &=1_{y'\ne y_0} J(y_0,y') +
[\lambda(x_0,y_0)-\lambda(y_0)]\delta_{y_0}(y')\\
&= J^{\pi_0}_y(y_0,y') +
[\lambda(x_0,y_0)-\lambda^{\pi_0}(x_0,y_0)] \delta_{y_0}(y')\\
\endaligned$$
for the transport cost $T_c( J^{\lambda}_x(x_0, \cdot),
J^{\lambda}_y(y_0, \cdot))$. Since $J^{\pi_0}((x_0,y_0),\cdot) +
[\lambda(x_0,y_0)-\lambda^{\pi_0}(x_0,y_0)] \delta_{(x_0,y_0)}$ is a
coupling of $J^{\lambda}_x(x_0, \cdot),  J^{\lambda}_y(y_0, \cdot)$,
we have
$$
T_c(J^{\lambda}_x(x_0, \cdot), J^{\lambda}_y(y_0, \cdot))\le
T_c(J^{\pi_0}_x(x_0,\cdot),
J^{\pi_0}_y(y_0,\cdot))+[\lambda(x_0,y_0)-\lambda^{\pi_0}(x_0,y_0)]c(x_0,y_0).
$$
Then for the Markov generator $\LL^{\pi,\lambda}$ with coupling
jumping  rates kernel $J^{\pi,\lambda}$, it is again a coupling
generator of $\LL$, and
$$\aligned
&\LL^{\pi,\lambda}c(x_0,y_0) \sum_{(x',y')\in S^2}
[c(x',y')-c(x_0,y_0)]J^{\pi,\lambda}((x_0,y_0),(x',y'))\\
&= T_c(J^{\lambda}_x(x_0, \cdot), J^{\lambda}_y(y_0, \cdot)) -\lambda(x_0,y_0) c(x_0,y_0)\\
&\le T_c(J^{\pi_0}_x(x_0,\cdot),
J^{\pi_0}_y(y_0,\cdot))+[\lambda(x_0,y_0)-\lambda^{\pi_0}(x_0,y_0)]c(x_0,y_0)
-\lambda(x_0,y_0) c(x_0,y_0)\\
 &=\LL^{\pi_0}c(x_0,y_0).
\endaligned$$
As $\LL^{\pi_0}$ is already an optimal coupling, the above
inequality must be equality. Thus $\LL^{\pi,\lambda}$ is also an
optimal coupling.

As $\lambda^{\pi_0}(x_0,y_0)\le \lambda(x_0) +\lambda(y_0)$, we
finish the proof of the theorem.
 \nprf

Following \cite{Oll09, LLY11} and Theorem \ref{thm22}, we can define
the Ricci curvature in term of generator:

\bdef\label{def_Ricci} {\it The Ricci curvature of $\LL$ w.r.t. a
metric $d$ is defined as
$$
{\rm Ric}_{(x,y)}(\LL,d)=- \inf_{\LL^\pi} \frac{\LL^\pi
d(x,y)}{d(x,y)}= \frac{[\lambda(x)+\lambda(y)]d(x,y)
-T_d(J(x,\cdot), J(y,\cdot))}{d(x,y)}
$$
for every $x\ne y$ in $S$, where the infimum is taken over all
coupling Markov generators of $\LL$, and $J(x,\cdot), J(y,\cdot)$
are jumping  rates measures with $J(x,x)=\lambda(y),
J(y,y)=\lambda(x)$.

The Ricci curvature at $(x,y)\in S^2\backslash \triangle$ of the graph $G$ is defined by
$$
Ric_{(x,y)}^G:= Ric_{(x,y)}(\Delta, d_G).
$$
}\ndef

We introduce now {\it length metric}.
 A positive function $w:E\rightarrow (0,+\infty)$
defined on the edge set $E$ is called {\it length function}, if
$w(x,y)=w(y,x)$ for any $e=(x,y)\in E$. Given the length function
$w$,  the {\it $w$-length} of a path $\gamma_{xy}$ from $x$ to $y$
is defined by
$$ |\gamma_{xy}|_w:=\sum_{e\in \gamma_{xy}}w(e).
$$ The length distance between $x,y$ associated with $w$  is
defined by
$$
d_w(x,y):=\min_{\gamma_{xy}} |\gamma_{xy}|_w.
$$
When $w\equiv 1$ on $E$, $d_w=d_G$ is the natural graph distance on
$S$.

\bcor\label{cor24} For the length metric $d_w$, $Ric(\LL, d_w)\ge
\kappa$ if and only if for any $x\sim y$, there is a coupling
measure $J^\pi((x,y),\cdot)$ of $J(x,\cdot)$ with
$J(x,x)=\lambda(y)$ and $J(y,\cdot)$ with $J(y,y)=\lambda(x)$, such
that
$$
\sum_{(x',y')\in S^2} J^\pi((x,y), (x',y'))[d_w(x',y')-d_w(x,y)] \le
-\kappa d_w(x,y).
$$
In other words, the best lower bound $\kappa$ of Ricci curvature of
$\LL$ w.r.t. $d_w$ can be identified as \beq\label{cor24a}\aligned
\kappa&=\inf_{x\sim y} \frac{[\lambda(x)+\lambda(y)]
[d_w(x,y)-W_{1,d_w}(\tilde J_x, \tilde J_y)] }{d_w(x,y)}\\
 \endaligned
\nneq where
$$\aligned
\tilde J_x(x') &= \frac {1_{x'\ne x}J(x,x') +
\lambda(y)\delta_x(x')}{\lambda(x)+\lambda(y)},\quad  \tilde J_y(y')
= \frac {1_{y'\ne y}J(y,y') +
\lambda(x)\delta_y(y')}{\lambda(x)+\lambda(y)}.\\
\endaligned$$
 \ncor

\bprf The necessity follows immediately from Theorem \ref{thm21}
and Theorem \ref{thm22}. For the sufficiency,  for any $x\ne y$ which are
not neighbors, let $\gamma_{x,y}=\{(x_k,x_{k+1});\ k=0,\cdots,n-1\}$
be a geodesic from $x$ to $y$ ($x_0=x,x_n=y$) in the length metric
$d_w$, i.e. $|\gamma_{xy}|_w=d_w(x,y)$. Take some coupling jumping
rate $\lambda(x,y)\ge \max\{\lambda(x_i)+\lambda(x_j); i\ne j\}$ and
consider the positive measures of the same mass $\lambda(x,y)$:
$$
J_k(x_k')=1_{x'_k\ne x_k}J(x_k,x'_k) + [\lambda(x,y)-\lambda(x_k)]
\delta_{x_k}(x'_k),\ k=0,\cdots, n
$$
and $\tilde J_k =J_k/\lambda(x,y)$, the corresponding normalized
probability measures. By Theorem \ref{thm22}, for an optimal
coupling generator $\LL^\pi$ in the sense that for every
$(x_0,y_0)\in S^2$, $\LL^\pi d_w(x_0,y_0)$ attains the minimum among
all coupling generators,
$$\aligned
\LL^\pi d_w(x,y) &= T_{d_w}(J_0, J_n) - \lambda(x,y) d_w(x,y)\\
&=\lambda(x,y) \left[W_{1,d_w}(\tilde J_0,\tilde J_n)
-\sum_{k=0}^{n-1}
d_w(x_k,x_{k+1})\right]\\
&\le \lambda(x,y)\sum_{k=0}^{n-1} \left[W_{1,d_w}(\tilde J_k,\tilde
J_{k+1}) -
d_w(x_k,x_{k+1})\right]\\
&=\sum_{k=0}^{n-1} \LL^\pi d_w(x_k,x_{k+1})\\
&\le -\kappa \sum_{k=0}^{n-1} d_w(x_k,x_{k+1})=-\kappa d_w(x,y),\\
\endaligned$$
where the last inequality follows by our condition on the neighbors
$(x_k,x_{k+1})$. That implies $Ric(\LL, d_w)\ge \kappa$ by Theorem
\ref{thm21}.

Finally (\ref{cor24a}) follows from the previous conclusion and the
construction of the optimal couplings in Theorem \ref{thm22}.
 \nprf

From Theorem \ref{thm22}, we derive immediately
 \bcor If $\lambda(x)=1$ for all $x\in S$ (i.e. $J$ is a probability transition
 kernel), then $Ric(\LL, d)\ge \kappa$ if and only if
 $$
2[W_{1,d}(\tilde J_x, \tilde J_y) -d(x,y)]\le -\kappa d(x,y),
 $$
where $\tilde J_x(x')=\frac 12 [1_{x'\ne x}J(x,x')+ \delta_x(x')]$,
in other words the $\alpha$-Ricci curvature lower bound
$\kappa^{(\alpha)}$ with $\alpha=1$ in the sense of Lin-Lu-Yau
\cite{LLY11} coincides with the lower bound of Ricci curvature
defined in (\ref{Ric}).
 \ncor

\brmk{\rm  If $d(x,y)=1_{x\ne y}$ is the discrete metric, then
$$W_{1,d}(\nu_1,\nu_2)=\|\nu_1-\nu_2\|_{tv}:=\sup_{A\subset
S}|\nu_1(A)-\nu_2(A)|.$$ A quite natural question is: in the
uniformly ergodic case, whether one has the exact exponential
convergence:
$$
\|P_t(x,\cdot)-P_t(y,\cdot)\|_{tv} \le e^{-\kappa t},\ t\ge0, x,y\in
S
$$
for some constant $\kappa>0$. That is $Ric(\LL,d)\ge \kappa>0$. The
answer is

\bcor\label{cor28} For the discrete metric $d(x,y)=1_{x\ne y}$ and
some positive constant $\kappa>0$, $Ric(\LL,d)\ge \kappa$ if and only if the
graph diameter $D_G=\sup_{x,y\in S} d_G(x,y) \le 2$, and
for any $x\ne y$,
$$
\sum_{x'\in S} J(x,x')\wedge J(y,x')\ge \kappa,
$$
where $J(x,x):=\lambda(y)\ge J(y,x)$ and $J(y,y):=\lambda(x)\ge
J(x,y)$. \ncor

\bprf When the graph diameter $D_G\ge 3$, for $(x,y)\in
S^2$ such that $d_G(x,y)\ge 3$, since $T_d(J(x,\cdot),
J(y,\cdot))=[\lambda(x)+\lambda(y)]$ when $J(x,x):=\lambda(y)$ and
$J(y,y):=\lambda(x)$,
 we get that the Ricci curvature ${\rm Ric}_{(x,y)}(\LL,d)$ w.r.t. the
discrete metric  equals to zero by Theorem \ref{thm22}.

Assume now $D_G(S)\le 2$. For any $x\ne y$,  since
$$T_d(J(x,\cdot), J(y,\cdot))=\sum_{x'\in A}(J(x,x')-J(y,x')),$$ where
$A=\{x': J(x,x')>J(y,x')\}$, we get by using the expression in
Definition \ref{def_Ricci},
$$
\aligned {\rm
Ric}_{(x,y)}(\LL,d)&=[\lambda(x)+\lambda(y)]-T_d(J(x,\cdot),
J(y,\cdot))\\
&=[\lambda(x)+\lambda(y)] - J(x,A) +J(y,A)\\
 &=J(x,A^c)+J(y,A)\\
&=\sum_{x'\in S} J(x,x')\wedge J(y,x'),
\endaligned$$
which is the desired result.
 \nprf
 }\nrmk

Notice that if $\LL$ is $\mu$-symmetric or equivalently $\mu(x)J(x,y)=\mu(y)J(y,x)$ for all $(x,y)\in S^2$, if $Ric(\LL,d)\ge \kappa$, the infimum $\lambda_1$ of the spectrum of $-\LL$ above zero on $L^2(S,\mu)$ (i.e. the spectral gap of $\LL$) satisfies (\cite{Wu04a})
\beq\label{713a}
\lambda_1\ge \kappa.
\nneq

\bexa {\bf (complete graph) }{\rm Let $(S,E)$ be a complete graph
with $N\ge 2$ vertices, i.e. $x\sim y$ for any two different
vertices $x,y$. Consider the Laplacian $\displaystyle \Delta f(x)=
\frac 1{N-1}\sum_{y\in S} (f(y)-f(x))$.

For this model the graph metric $d_G(x,y)$ is the discrete metric $1_{x\ne y}$. By
Corollary \ref{cor28},
$$
Ric_{(x, y)}^G=Ric_{(x,y)}(\Delta, d_G)=\frac N{N-1},\ (x,y)\in S^2\backslash \triangle.
$$
Remark that the first eigenvalue $\lambda_1$ of $-\Delta$ above zero
is $\frac {N}{N-1}$.

When $N=2$, we get the affirmation for the graph of two vertices in
the Introduction.
 } \nexa

\bexa {\bf (star-graph)} {\rm $S=\{o,x_1,\cdot,x_N\}$ $(N\ge 2)$
where $o\sim x_k$ but the only neighbor of $x_k$ is $o$,
$k=1,\cdots, N$. For the Laplacian $\Delta$, the spectral gap $\lambda_1=1$.

By Corollary \ref{cor24}, the best lower bound $\kappa_G$ of Ricci
curvature for the graph metric $d_G$ is
$$
\kappa_G =2-T_{d_G}(J(o,\cdot),J(x_1,\cdot)),
$$
where $J(o,o)=J(x_1,x_1)=1$. By transporting the mass of
$x_2,\cdots,x_N$ of $J(o,\cdot)$ to $x_1$, we see that
$T_{d_G}(J(o,\cdot),J(x_1,\cdot))\le 2\times \frac{N-1}{N}$. On the
other hand taking $f(x_1)=-1$, $f(o)=0$ and $f(x_k)=1$ for $k\ge 2$
we have by Kantorovitch's duality characterization
$$
T_{d_G}(J(o,\cdot),J(x_1,\cdot))\ge Jf(o)-Jf(x_1)=\frac{N-1}{N}
-\frac 1N - (-1)=\frac{2(N-1)}{N}.
$$
The transport above is therefore optimal and then
$$
\kappa_G=\frac{2}{N}
$$
which is far smaller than the spectral gap $\lambda_1=1$ for big $N$.

However
we can construct a new metric $d$ for which $Ric(\Delta, d)\ge
\lambda_1=1$. This new metric is
$$
d(x_i,x_j)=2, \ i\ne j;\ d(o,x_j)=\frac{2(N-1)}{N}\ge 1.
$$
Using the coupling generator $\LL^\pi$ of $\Delta$ with jumping
rates
$$
J^\pi((x_i,x_j), (o,o))=1 \ (i\ne j), J^{\pi}((o,x_i), (o,o))=1,
J^\pi((o,x_i), (x_k,x_i))=\frac 1N, \ \forall i,k
$$
and symmetrically for $J^\pi((x_i,o),\cdot)$, we have
$$
\LL^\pi d(x,y) = - d(x,y), \ \forall (x,y)\in S^2
$$
which means that $Ric(\Delta, d)\ge 1=\lambda_1$, the optimal Ricci curvature lower bound.
 } \nexa

 \subsection{Myer's diameter theorem}
The following result is the continuous time counterpart of
Ollivier's theorem \cite{Oll09} in the discrete time case:

 \bprop\label{prop27} If $Ric(\LL,d)\ge \kappa>0$, then
for all $x\ne y$ in $S$, \beq\label{prop27a} d(x,y) \le \frac
1\kappa\left(\sum_{x'\sim x} J(x,x') d(x,x') + \sum_{y'\sim y}
J(y,y') d(y,y')\right). \nneq In particular if $Ric(\LL,d_G)\ge
\kappa>0$ and $M=\sup_x \lambda(x)<+\infty$, then the diameter
$D_G=\sup_{x,y\in S}d_G(x,y)$ satisfies \beq\label{prop27b} D_G\le
\frac{2M}{\kappa}. \nneq
 \nprop
The estimate (\ref{prop27b}) is sharp as seen for the two vertices
graph with $\Delta$.

\bprf By Theorems \ref{thm21} and \ref{thm22} (by taking
$\lambda(x,y)=\lambda(x)+\lambda(y)$),
$$
\kappa d(x,y) \le [\lambda(x)+\lambda(y)]d(x,y) - T_{d}(J(x,\cdot),
J(y,\cdot))
$$
where $J(x,x)=\lambda(y)$ and $J(y,y)=\lambda(x)$ (as specified in
(\ref{thm22a})). By the triangular inequality,
$$\aligned
&[\lambda(x)+\lambda(y)]d(x,y)=T_d([\lambda(x)+\lambda(y)]\delta_x,
[\lambda(x)+\lambda(y)]\delta_y)\\
&\le T_d( [\lambda(x)+\lambda(y)]\delta_x, J(x,\cdot)) + T_d(
[\lambda(x)+\lambda(y)]\delta_y, J(y,\cdot))+
T_d(J(x,\cdot), J(y,\cdot))\\
&= \sum_{x'\sim x} J(x,x')d(x,x') + \sum_{y'\sim y} J(y,y')d(y,y') +
T_d(J(x,\cdot), J(y,\cdot))
\endaligned
$$
Plugging it into the previous inequality we obtain (\ref{prop27a}). Finally (\ref{prop27b}) is a direct
consequence of  (\ref{prop27a}). \nprf

\bexa{\rm Let $S=\{k\in \nn;\ k\le N\}$ ($3\le N\in\nn$) be equipped
with graph metric $d_G(x,y)=|y-x|$. Consider the generator
$$
\LL f(n)=b [f(n+1)-f(n)] + n [f(n-1)-f(n)], n\in S
$$
where $b>0$, $f(-1):=f(0), f(N+1):=f(N)$. By Corollary \ref{cor2-10}
below, $Ric(\LL, d_G)\ge 1$. By (\ref{prop27a}) in the Myer's type
theorem above with $x=0, y=N$,
$$
D_G\le N+b.
$$
This upper bound of the diameter becomes sharp by letting $b\to0+$
(as $D_G=N$). That shows again the sharpness of Proposition
\ref{prop27}.
 }\nexa
\subsection{Superposition and tensorization}
\bprop If the Ricci curvature of Markov generators $\LL_i,
i=1,\cdots, N$ of nearest-neighbor type (as $\LL$) w.r.t. some
metric $d$ on $S$ are bounded from below by $\kappa$, then for
positive constants $\alpha_1,\cdots, \alpha_N$,
$$
Ric\left(\sum_{i=1}^N \alpha_i \LL_i, d\right)\ge \kappa
\sum_{i=1}^N \alpha_i.
$$
\nprop

\bprf By Theorem \ref{thm22}, there are coupling Markov generators
$\LL_1^\pi,\cdots, \LL^\pi_N$ of $\LL_1,\cdots,\LL_N$ such that
$$
\LL^\pi_i d(x,y)\le -\kappa  d(x,y), \ \forall (x,y)\in S^2,\ i=1,2,\cdots,N.
$$
As $\sum_{i=1}^N \alpha_i \LL_i^\pi$ is a coupling Markov generator
of $\sum_{i=1}^N \alpha_i \LL_i$ and
$$
(\sum_{i=1}^N \alpha_i \LL_i^\pi )d(x,y) \le -\kappa [\sum_{i=1}^N
\alpha_i] d(x,y),
$$
where the desired result follows by Theorem \ref{thm21}.
 \nprf

Now we show that the Ricci curvature lower bound is dimension-free.
Let $(S_1,E_1),
\\
\cdots,(S_N,E_N)$ be $N$ graphs and consider the product graph
$S=\prod_{i=1}^N S_i$: two vertices $x=(x_1,\cdots, x_N)$ and
$y=(y_1,\cdots,y_N)$ are adjacent if and only if $\exists j$ such
that $x_j\sim y_j$ and $x_i=y_i$ for all $i\ne j$.

\bprop\label{prop-tensorization} Let $\LL_i$ be a Markov generator
on $S_i$ with jumping rates kernel $J_i$ of nearest-neighbor type
where $i=1,\cdots, N$, and consider the Markov generator on the
product graph $S=\prod_{i=1}^N S_i$:
\beq\label{prop-tens1} \LL
f(x)= \sum_{i=1}^N \sum_{y_i \in S}[f(x^{y_i})-f(x)] J_i(x_i,y_i),\
x=(x_1,\cdots, x_N)\in S,\nneq
where $(x^{y_i})_j=x_j$ for all $j\ne
i$ and $(x^{y_i})_i=y_i$.
If $Ric(\LL_i,d_i)\ge \kappa$ for all $i=1,\cdots,N$ where
$d_i$ is some metric on $S_i$, then $Ric(\LL, d_{L^1})\ge \kappa$
where \beq\label{L1metric} d_{L^1}(x,y)=\sum_{i=1}^N d_i(x_i,y_i)
\nneq is the $L^1$-metric on $S=\prod_{i=1}^N S_i$. \nprop

The generator $\LL$ in (\ref{prop-tens1}) is often denoted by
$\oplus_{i=1}^N \LL_i$. The semigroup $(P_t)$ generated by $\LL$ is
the tensorization of the semigroups $(P_t^{(i)})$ generated by
$\LL_i$, i.e.
$$
P_t(x,y) =\prod_{i=1}^N P_t^{(i)}(x_i,y_i).
$$

\bprf By Theorem \ref{thm21}, for every $i$ there is a coupling
Markov generator $\LL_i^\pi$ of $\LL_i$ such that
$$
\LL_i^\pi d_i(x_i,y_i)\le -\kappa d_i(x_i,y_i), \ (x_i,y_i)\in
S_i^2.
$$
Letting $\LL^\pi=\oplus_{i=1}^n \LL_i^\pi$ which is a coupling
Markov generator of $\LL$, we have
$$
\LL^\pi d_{L^1}(x,y) =\sum_{i=1}^N \LL_i^\pi d_i(x_i,y_i) \le
-\kappa \sum_{i=1}^N  d_i(x_i,y_i),
$$
which implies the desired result by Theorem \ref{thm21} again.
 \nprf

\subsection{Order-preserving coupling}
Let $S$ be equipped with a partial order $\preceq$ such that neighbors are comparable: if $x\sim y$ then either $x\preceq y$ or
$y\preceq x$.

For two probability measures $\nu_1,\nu_2$ on $S$, we say that
$\nu_1\preceq\nu_2$, if there are two random variables $X_1$ and
$X_2$ valued in $S$ with laws $\nu_1$ and $\nu_2$, such that $X_1\preceq
X_2, a.s. $. In that case we say that $(X_1,X_2)$ is an ordering
coupling of $\nu_1$ and $\nu_2$.

For two positive measures $\nu_1, \nu_2$ of the same mass $m$, we
say that $\nu_1\preceq\nu_2$, if $\nu_1/m \preceq \nu_2/m$.

\blem\label{lem25} Assume that $d_w$ is a length metric on $S$ such
that there is an increasing function $h$ w.r.t. the order $\preceq$
so that \beq\label{lem25a} d_w(x,y)= h(y) - h(x), \ \forall x\preceq
y. \nneq Given two probability measures $\nu_1\preceq \nu_2$, every
ordering coupling of $\nu_1,\nu_2$ is optimal for
$W_{1,d_w}(\nu_1,\nu_2)$ and
 \beq\label{lem25b}
W_{1,d_w}(\nu_1,\nu_2)=\sum_{x\in S} h(x)[\nu_2(x)-\nu_1(x)].\nneq
 \nlem

\bprf Let $X_1\preceq X_2$ be two random variables valued in $S$,
with laws $\nu_1$ and $\nu_2$. We have
$$
W_{1,d_w}(\nu_1,\nu_2) \le \ee d_w(X_1,X_2) = \ee [h(X_2)-h(X_1)]
=\sum_{x\in S} h(x)[\nu_2(x)-\nu_1(x)].
$$
The converse inequality holds by Kantorovitch duality  for
$\|h\|_{Lip(d_w)}=1$, where $ \|h\|_{Lip(d_w)}:=\sup_{x,y\in S}\frac{|h(x)-h(y)|}{ d_w(x,y)}$.
\nprf

\bprop\label{prop28} Assume that $(S,E)$ be an ordered graph
equipped with the length metric $d_w$ satisfying (\ref{lem25a}). If
for every $x\sim y$ with $x\preceq y$,
$$
J(x,\cdot) \preceq J(y,\cdot) \text{ with }\ J(x,x)=\lambda(y),
J(y,y)=\lambda(x),
$$
then $Ric(\LL, d_w)\ge \kappa$ if and only if for every $x\sim y$
with $x\preceq y$,
$$
\LL h(y) - \LL h(x)\le -\kappa [h(y)-h(x)].
$$
\bprf By (\ref{cor24a}) in Corollary \ref{cor24}, $Ric(\LL, d_w)\ge \kappa$ iff (if and only if)
for every $x\sim y$ with $x\preceq y$,
$$
[\lambda(x) + \lambda (y)] (W_{1,d_w}(\tilde J_x, \tilde J_y) -
d_w(x,y))\le -\kappa d_w(x,y) = -\kappa (h(y)-h(x)),
$$
where $\tilde J_x=J(x,\cdot)/(\lambda(x) +\lambda(y))$ is the
corresponding normalized probability measure. As $\tilde J_x \preceq
\tilde J_y$, we have by Lemma \ref{lem25}
$$\aligned
&[\lambda(x) + \lambda (y)] (W_{1,d_w}(\tilde J_x, \tilde J_y) -
d_w(x,y))\\
&= [\lambda(x) + \lambda (y)] \left\{\sum_{z} h(z) (\tilde J_y(z)-
\tilde J_x(z)) - (h(y) - h(x))\right\}\\
&=\LL h(y) - \LL h(x), \endaligned$$ where the desired result
follows.
 \nprf
 \nprop

Now we turn to the overworking model: birth-death processes.

\bexa{\rm (Birth-death processes) Let $S=\nn\bigcap [0,D]$ where
$D\in \nn^*\cup \{+\infty\}$, $x\sim y$ iff $|y-x|=1$. Its graph
distance coincides with the Euclidian one. Consider the birth-death
process \beq\label{bd1} \LL f(n)=b_n (f(n+1)-f(n)) + a_n
(f(n-1)-f(n)), \ n\in S, \nneq where $a_0=0, a_n>0 (n\ge1)$ and
$b_n>0 (0\le n\le D-1), b_D=0$ (if $D$ is finite).

\bcor\label{cor2-10} For the birth-death generator $\LL$ given in
(\ref{bd1}), $Ric(\LL, d_G)\ge \kappa$ if and only if
\beq\label{bd2} (a_{n+1}-a_n)-(b_{n+1}-b_n)\ge \kappa, \ \forall
n\in [0,D-1]\cap \nn.\nneq \ncor

\bprf Let $h(n)=n, \ n\in [0,D]\cap\nn$. For this model
$J(n,n+1)=b_n, J(n,n-1)=a_n$. Letting
$J(n+1,n+1)=\lambda(n)=a_n+b_n$ and $J(n,n)=\lambda(n+1)$, we have
$J(n,\cdot)\preceq J(n+1,\cdot)$. By Proposition \ref{prop28},
$Ric(\LL, d_G)\ge \kappa$ if and only if
$$
\LL h(n+1) -\LL h(n)= (b_{n+1} - a_{n+1}) - (b_n-a_n)\le -\kappa, \
\forall n\in [0,D-1]\cap\nn,
$$
which is the desired result.
\nprf
 }\nexa

The condition (\ref{bd2}) with $\kappa>0$ was introduced by Caputo, Dai Pra and Posta \cite{CdPP09} as the counterpart of Bakry-Emery's positive curvature condition. They established the exponential convergence in entropy of the birth-death process under (\ref{bd2}) together with the non-decreasingness of $(a_n)$ and the non-increasingness of $(b_n)$.

\bexa {\bf ($M/M/\infty$ queue and Poisson measure)} {\rm Consider
the $M/M/\infty$ queue: the birth-death process valued in $S=\nn$
with the generator
$$
\LL f(n)=\lambda[f(n+1)-f(n)] + n[f(n-1)-f(n)],\ n\in\nn,
$$
whose unique invariant probability measure is the Poisson
distribution $\mu$ with parameter $\lambda>0$. This example is the
counterpart of the Ornstein-Uhlenbeck process, in  jumping Markov
processes.

By Corollary \ref{cor2-10},
$$
Ric(\LL, d_G)\ge 1,
$$
which is the exact Ricci curvature lower bound for $\lambda_1=1$.
This fact is well known since $D_+P_tf=e^{-t}P_t D_+f$, where
$D_+f(n):=f(n+1)-f(n)$. }\nexa

\bexa {\bf (Discrete cube)} {\rm Consider the product measure
$\mu^{\otimes n}$ on the product graph $S =\{0,1\}^n$ equipped with
the graph metric $d_{G} (x,y)= \sum_{k=1}^n 1_{x_i\ne y_i}$ (the
Hamming metric), where $\mu (1)=p$ and $\mu(0)=1-p=:q$ ($0<p<1$).
Let $J_m$ be a jumping rate kernel on $\{0,1\}$ satisfying the
detailed balance condition: $\mu(0)J_m(0,1)=\mu(1)J_m(1,0)$.
Consider the corresponding generator on $\{0,1\}$:
$$
\begin{aligned}
\LL_m f(x) = 1_{\{0\}}(x) [f(1) - f(0)] J_m(0,1) + 1_{\{1\}}(x)
[f(0)- f(1)] J_m(1,0).
\end{aligned}
$$
By Corollary \ref{cor28}, $Ric_{(0,1)}(\LL_m, d)=J_m(0,1)+J_m(1,0)$, where
$d$ is the graph metric on $\{0,1\}$ (i.e. $d(0,1)=1$). On the other
hand, by direct calculus the spectral gap $\lambda_1$ of $\LL_m$
equals to $J_m(0,1)+J_m(1,0)$, too.

Consider the generator $\LL = \oplus_{i=1}^n \LL_i$ on $S=\{0,1\}^n$,
where $\LL_i=\LL_m$. By the tensorization (Proposition
\ref{prop-tensorization}),
$$
Ric(\LL,d_G)\ge J_m(0,1)+J_m(1,0).
$$

}\nexa

\bexa {\bf (Binomial distribution)} {\rm Let
$S=I_n:=\{0,1,\cdots,n\}$ equipped with the graph metric as the
Euclidean metric and $\mu$ the binomial law $\BB(n,p)$ on
$S=\{0,1,\cdots,n\}$, where $0<p<1$. Consider the generator
$$
\LL_{bin} g (y)=p(n-y)[g(y+1)-g(y)]+(1-p)y[g(y-1)-g(y)], \ y\in S
$$
which is symmetric on $L^2(\mu)$. By Corollary \ref{cor2-10},
$$
Ric(\LL_{bin}, d_G)\ge 1.
$$
Further since $\LL_{bin} h=-h$ where $h(y)=y-\mu(h)$ is increasing,
we see that $\lambda_1=1$. In other words the estimate of Ricci
curvature above is sharp.

}

\nexa

 \bexa{\bf (Ferromagnetic spin model)}
{\rm Let $S=\{-1,1\}^N$, equipped with the graph metric
$d_G(x,y)=\sum_{k=1}^N 1_{x_k\ne y_k}$ (the Hamming metric).
Consider
$$
\mu(x) = \frac {e^{-V(x)}}{\sum_{y\in S} e^{-V(y)}}
$$
where $V(x)= \sum_{i<j} \beta_{ij} x_i x_j $ with $\beta_{ij}\le 0$.
The Gibbs algorithm for sampling the high dimensional distribution
$\mu$ is given by
$$
\LL f(x)= \sum_{i=1}^N (\bar \mu_i(f)(x)-f(x))
$$
where $\bar \mu_i(x'|x)=\prod_{j\ne i}\delta_{x_j}(x_j') \cdot
\mu_i(x_i'|x)$, and
$$
\mu_i(x'_i|x)= \frac{e^{-V(x^{x_i'})}}{e^{-V(x^{i+})} + e^{-V(x^{i-})}},\
$$
is the conditional distribution of $x_i$ knowing $(x_j)_{j\ne i}$
under $\mu$. Here $ x^{i\pm}_j=x_j (j\ne i), x^{i\pm}_i=\pm 1$.
Since
$$\aligned
c_{ij}(x) :&=\mu_i(1|x^{j+}) -\mu_i(1|x^{j-}) \\
&= \frac{e^{2\sum_{k\ne i,j} \beta_{ik} x_k } (e^{ -2\beta_{ij}}-
e^{2\beta_{ij}})}{ (1+ e^{2\sum_{k\ne i,j} \beta_{ik} x_k
-2\beta_{ij}}) (1+ e^{2\sum_{k\ne i,j} \beta_{ik} x_k +
2\beta_{ij}}) }\ge 0, \endaligned$$  we can apply Proposition \ref{prop28} for
$d_w=d_G$ and
$$
h(x)=\frac 12 \sum_{i=1}^N x_i.
$$
Then $Ric(\LL, d_G)\ge \kappa$ if and only if for any $j=1,\cdots,N$
and $x\in S$ with $x_j=-1$,
$$
\LL h(x^{j+}) -\LL h(x) \le -\kappa. $$ But noting that
$\mu_j(h|x^{j+})=\mu_j(h|x)$, we have
$$\aligned\LL h(x^{j+}) -\LL h(x) & = \sum_{i:i\ne j}
\{[\mu_i(h|x^{j+}) -h(x^{j+}) ] - [\mu_i(h|x)-h(x)]\} -(h(x^{j+})-h(x))
\\
&=\frac 12\sum_{i: i\ne j}[(\mu_i(1|x^{j+})-\mu_i(-1|x^{j+})-x_i) - (\mu_i(1|x)-\mu_i(-1|x)-x_i)]-1 \\
&= \sum_{i:i\ne j} c_{ij}(x) -1,
\endaligned
$$
where the last equality follows from $$\mu_i(-1|x)- \mu_i(-1|x^{j+})= \mu_i(1|x^{j+})- \mu_i(1|x)=c_{ij}(x).$$
Therefore $Ric(\LL,d_G)\ge \kappa$ if and only if \beq \sup_{x\in
S}\max _{j}\sum_{i: i\ne j} c_{i,j}(x)\le 1-\kappa. \nneq If
$\kappa>0$, this becomes the famous Dobrushin's uniqueness condition
for Gibbs measure if one takes $\sup_{x}$ inside the sum above.

The sufficiency of the Dobrushin's uniqueness condition for the
exponential convergence in $W_1$-metric was found by the third
named author in \cite{Wu06} and Ollivier \cite{Oll09} (from the viewpoint of Ricci curvature).
 }
 \nexa

\section{Comparison with a birth-death process}

Often the Ricci curvature of $\LL$ w.r.t. the graph metric $d_G$ is
bounded from below by zero or a negative number.  The purpose of
this and the next two sections is to introduce a comparison technique to find a metric
$d$ so that $Ric(\LL, d)\ge \kappa>0$ even when $Ric_{(x,y)}(\LL, d_G)\le 0$.

\subsection{Ricci curvature of birth-death processes w.r.t. a general length metric}

For a general birth-death process valued in $S=[0,D]\cap\nn$,
\beq\label{bd-generator} \LL f(n)= a_n(f(n-1)-f(n)) + b_n
(f(n+1)-f(n)), \ n\in [0,D]\cap\nn,\nneq where $a_n>0\ (n\ge1),
b_n>0\ (n\le D-1)$ are respectively the death and birth rates ($-1$
is identified with $0$ and $D+1$ is identified as $D$ if $D$ is
finite). In the following $a_0:=0$ and $b_D:=0$ if $D$ is finite.

We always assume that the Markov process generated by $\LL$ is
conservative. This process is always symmetric w.r.t. the invariant
probability measure
$$
\mu_0=\frac 1C,\ \mu(n)=\frac 1C\frac{b_0\cdots b_{n-1}}{a_1\cdots
a_n},\ 1\le n\le D,
$$
where $C>0$ is the normalization constant. Any length metric $d_w$
is determined by an increasing function $h$ on $S$:
$w(n,n+1)=h(n+1)-h(n)=:D_+h(n)$. Since
$$\aligned
\LL& h(n+1)-\LL h(n)= b_{n+1} D_+h(n+1) - a_{n+1} D_+h(n) - b_n
D_+h(n)+
a_n D_+h(n-1)\\
&=-[(a_{n+1} D_+h(n)- a_n D_+h(n-1))-( b_{n+1} D_+h(n+1)-  b_n
D_+h(n))]
\endaligned
$$
by Proposition \ref{prop28}, we get immediately

\bprop\label{prop-bd} For the birth-death process with generator
$\LL$ above and for the length metric $d_w$ determined by an
increasing function $h$, we have $Ric(\LL, d_w)\ge \kappa$ if and
only if for all $n\in S\backslash \{D\}$, \beq\label{bd3} (a_{n+1}
D_+h(n)- a_n D_+h(n-1))-( b_{n+1} D_+h(n+1)- b_n D_+h(n))\ge \kappa
D_+h(n) \nneq where $a_0:=0, b_D:=0$ if $D$ is finite. \nprop

\brmk{\rm If we consider $a_n'=a_nD_+h(n-1)=a_n w(n-1,n), \ b'_n=b_n
D_+h(n)=b_n w(n,n+1)$, the l.h.s. of (\ref{bd3}) is exactly
$$
(a_{n+1}'-a_n')-(b'_{n+1}-b'_n)
$$
the quantity appeared in Corollary \ref{cor2-10}. Notice that this
change of birth and death rates does not change the invariant
measure of $\LL$.
  }\nrmk

The proposition above is interesting (i.e. $\kappa> 0$) only if $g=-\LL h$ is
increasing and if such a function $h$ could be found. The following
proposition, based on a result of Liu and Ma \cite{LM09}, provides a
method to find such $h$.

\bprop\label{prop-bd2} Let $g:S\to \rr$ be an increasing
$\mu$-integrable function such that $\mu(g)=\sum_{k\in S}\mu(k)
g(k)=0$, and $h:S\to\rr$ be the solution to the Poisson equation
$-\LL h=g$, determined up to difference of a constant (\cite{LM09})
by \beq\label{Poisson-solution}
 D_{+}h(n-1)=h(n)-h(n-1) =\frac{\sum_{k\ge n}\mu(k)g(k)}{a_n\mu(n)},\ \forall
 1\le n\le D.
\nneq If \beq\label{Lip-norm} \|h\|_{Lip(g)}:=\sup_{1\le n\le D}
\frac{D_+h(n-1)}{D_+g(n-1)}=: K(g)<+\infty,\nneq then $Ric(\LL,
d_w)\ge \frac 1{K(g)}$, where $w(n,n+1)=h(n+1)-h(n)$. \nprop

\bprf Since
$$\sum_{k\ge n}\mu(k)g(k)=\sum_{i,j: i<n\le j\le D} (g(j)-g(i))
\mu(i)\mu(j)>0,
$$
$D_+h(n-1)>0$, i.e. $h$ is increasing. Moreover
$$
-\LL h(n+1) + \LL h(n) = g(n+1)-g(n)\ge \frac 1{K(g)} (h(n+1)-h(n)),
$$
where the desired result follows by Proposition \ref{prop-bd}.
 \nprf

 \brmk{\rm If $\LL$ admits an eigenfunction $h$
associated with the spectral gap $\lambda_1>0$, $h$ can be chosen as
increasing for which (\ref{bd3}) becomes equality and then the Ricci
curvature  ${\rm Ric}_{(m,n)}(\LL,d_w)=\lambda_1$ for the length
metric $d_w$ associated to $h$. The problem is that such
eigenfunction may not exist or may be very difficult to find
especially when $\lambda_1$ is unknown.

Under some extra conditions on the uniqueness of $\LL$ and the
integrability of $g$, Liu and Ma \cite{LM09} showed  that $K(g)$ is
exactly the norm of the Poisson operator $(-\LL)^{-1}$ in the space
of $g$-Lipschitzian functions with $\mu$-mean zero. Moreover
Chen's variational formula for the spectral gap $\lambda_1$
(\cite{Chen05}) says exactly
$$
\lambda_1 = \sup_{g\ \text{ increasing}} \frac{1}{K(g)}.
$$
In other words the Ricci curvature lower bound furnished by
Proposition \ref{prop-bd2} can attain or approach $\lambda_1$. }
\nrmk

Now we furnish a criterion for the exponential convergence of $(P_t)$
in $W_{1,d_G}$:
$$
W_{1,d_G}(P_t(x,\cdot), P_t(y,\cdot)) \le K e^{-\delta t} |x-y|,
t\ge0, x,y\in S.
$$

 \bcor\label{cor35} Assume that $m_\mu=\sum_{n\in S} n\mu(n)$ is
finite and for $g(x)=x-m_\mu$,
 \beq\label{cor35a} K(g):=\sup_{n\in S} \frac{\sum_{k\ge n}\mu(k)[k-m_\mu]}{a_n\mu(n)}<+\infty.  \nneq
\benu[(a)]

\item If
$$k(g)=\inf_{n\in S} \frac{\sum_{k\ge
n}\mu(k)[k-m_\mu]}{a_n\mu(n)}>0$$ then $$ W_{1,d_G}(P_t(x,\cdot),
P_t(y,\cdot)) \le \frac{K(g)}{k(g)} \exp\left(-\frac t{K(g)}\right)
\cdot |x-y|, \ t\ge0, x,y\in S.
$$
\item Assume that $k(g)=0$.
If moreover the Ricci curvature of $\LL$ w.r.t. the Euclidean metric
$d_G$ is bounded from below, i.e. for some constant $M\ge0$,
\beq\label{cor35b} (a_{n+1}-a_n)-(b_{n+1}-b_n)\ge -M,\ n\in S, \nneq
then for every $\alpha\in (0,1/M)$, for all $t\ge0,\ x,y\in S$, \beq
\label{cor35c} W_{1,d_G}(P_t(x,\cdot), P_t(y,\cdot)) \le
\frac{K(g)+\alpha}{\alpha}\exp\left(- \frac{1-\alpha
M}{K(g)+\alpha}\cdot t\right) d_G(x,y). \nneq \nenu
 \ncor

\bprf Let $h:S\to \rr$ be the function determined by $h(0)=0$ and
$$
D_+h(n-1):=h(n)-h(n-1)=\frac{\sum_{k\ge n}\mu(k)[k-m_\mu]}{a_n\mu(n)}, \
1\le n\le D_G.
$$
Then $\LL h(n)=-g(n),\ n\in S.$

{\bf (a).} By Proposition \ref{prop-bd2}, $Ric(\LL, d_w)\ge \frac
1{K(g)}$ where $w(n,n+1)=D_+h(n)$. Since $k(g)\le D_+h(n)\le K(g)$
for all $1\le n\le D-1$, we conclude easily the explicit exponential
convergence in the statement.

{\bf (b).} By the extra condition (\ref{cor35b}) we have
$$
D_+\LL g(n) = (b_{n+1}-a_{n+1}) -(b_{n}-a_{n}) \le M.
$$
Hence for any $\alpha\in(0,1/M)$, we have for any $n\le D-1$,
$$\aligned
D_+\LL (h+\alpha g)(n)&\le -1+\alpha M\\
&\le - \frac{1-\alpha M}{\sup_{k\le D_G-1} D_+(h+\alpha g)(k)}
D_+(h+\alpha g)(n)\\
&=- \frac{1-\alpha M}{K(g)+\alpha } D_+(h+\alpha g)(n).
\endaligned$$
Therefore by Proposition \ref{prop-bd}, for the metric $d_w$
associated with $w(n,n+1):=D_+h(n)+\alpha$,
$$Ric(\LL, d_w)\ge \delta:=
\frac{1-\alpha M}{K(g)+\alpha}.$$ As $\alpha\le w(n,n+1)\le
K(g)+\alpha$, we have
$$\aligned
W_{1,d_G}(P_t(x,\cdot), P_t(y,\cdot))&\le \frac{1}{\alpha}
W_{1,d_w}(P_t(x,\cdot), P_t(y,\cdot))\\
&\le \frac{1}{\alpha}e^{-\delta t} d_w(x,y)\\
&\le \frac{K(g)+\alpha}{\alpha}e^{-\delta t} d_G(x,y),
\endaligned$$
which is the desired result.
 \nprf

We now present several examples of birth-death processes which will
be served as reference models.

 \bexa\label{random walk}{\bf (Random walk on $[0,D]\cap \nn$ and uniform measure) }{\rm Consider the
 Laplacian $\Delta$ on $[0,D]\cap \nn$ where $D\in\nn^*$, i.e. $a_n=b_n=\frac 12$ for
 $0<n<D$ and $b_0=1, a_D=1$. Let
 $$
h(k) = -\cos \frac{k\pi}{D},
 $$
we have
$$
\Delta h(k)= - (1-\cos \frac \pi D) h(k),
$$
then by Proposition \ref{prop-bd}, for the metric $d_w$ with
$w(k,k+1)=h(k+1)-h(k)$, $Ric(\Delta, d_w)\ge 1-\cos \frac \pi D$
which is exactly the spectral gap $\lambda_1$ of $\Delta$.

 }\nexa

\bexa {\bf (Geometric measure I) }{\rm  A typical example is the
birth-death process
$$
\LL f(n) = b(f(n+1)-f(n)) +a (f(n-1)-f(n)), \ n\in \nn
$$
where $a>b>0$ and $f(-1):=f(0)$. Its unique invariant measure (up to
a constant factor) is $\mu(k)=(b/a)^k, k\in\nn$: the geometric
measure. Obviously w.r.t. the graph metric $d_G$, ${\rm Ric}_{(n,
n+k)}(\LL,d_G)=0$ for every $n,k\in \nn^*$ by applying Theorem
\ref{thm22} or Corollary \ref{cor2-10}. But its spectral gap
$\lambda_1$ is known: $\lambda_1= (\sqrt{a}-\sqrt{b})^2>0$
(\cite[Examples 9.22]{Chen92}). The problem is to find a metric $d$
such that $Ric(\LL, d)\ge \kappa>0$. In fact letting
$$
h(n)=\left(\sqrt{\frac{a}{b}}\right)^n,
$$
one has $\LL h(n)=-\lambda_1 h(n)$ for all $n\ge 1$ and $\LL
h(0)=h(0) [\sqrt{ab}-b]\ge -\lambda_1 h(0)$. Then
$$
\LL h(n+1) -\LL h(n) \le - \lambda_1 [h(n+1)-h(n)], \ n\in\nn,
$$
which implies by Proposition \ref{prop28} that $Ric(\LL, d_w)\ge
\lambda_1$, the best possible lower bound,  where
$w(n,n+1)=h(n+1)-h(n).$ } \nexa

Let $\tau_0:=\inf\{t\ge0; X_t=0\}$, the first hitting time to $0$.
Since $e^{\lambda_1 (t\wedge \tau_0)} h(X_{t\wedge\tau_0})$ is a
local martingale, then a supermartingale, we have
$$
\ee_n \exp\left((\sqrt{a}-\sqrt{b})^2 \tau_0\right)\le
\left(\sqrt{\frac{a}{b}}\right)^n, \ n\ge1.
$$
However this process is not exponentially convergent in the
Wasserstein metric $W_{1,d_G}$, i.e. there are no constants
$\kappa>0,\ K\ge 1$ such that
$$
W_{1,d_G}(P_t(x,\cdot), P_t(y,\cdot))\le K e^{-\kappa t} d_G(x,y), \
x,y\in\nn.
$$
In fact letting $D_+f(n)=f(n+1)-f(n)$, we have $D_+ \LL f=\LL
D_+f$, which implies $D_+P_tf=P_t D_+ f$. The exponential
convergence in $W_{1,d_G}$, if true, would imply that $P_t$ is
uniformly ergodic  (i.e. spectral gap exists in $L^\infty$). But
this process is not uniformly ergodic as well known.

\bexa {\bf (Geometric measure II)} {\rm Consider the birth-death
process
$$
\LL f(n) = p(n+1)(f(n+1)-f(n)) + n (f(n-1)-f(n)), \ n\in \nn,
$$
where $p\in(0,1)$. This is a model about the evolution of the number
of clients in a service center of infinite service capacity (such as
a great web site), such that new clients come in a rate proportional
to the number of clients in the center.

Its invariant measure is again the geometric measure
$\mu(n)=p^n(1-p)$. By Corollary \ref{cor2-10},
$$
Ric(\LL, d_G)\ge 1-p.
$$
Moreover for the increasing function $h(n)=n-\frac{p}{1-p}$ (note
that $\mu(h)=0$), since $\LL h(n)=-(1-p)h(n)$, then $\lambda_1=1-p$.
This shows the estimate above about the Ricci curvature is sharp.

 }\nexa

\subsection{Comparison with a generalized  birth-death process}

For a Markov process valued in the graph $(S,E)$ with generator
$\LL$, when its Ricci curvature lower bound  is not positive w.r.t.
the graph metric, we want to find some coupling generator $\LL^\pi$
and some increasing function $h_0:\nn\to \rr$ with $h_0(0)=0$ and
some positive constant $\kappa>0$ such that \beq \LL^\pi (h_0\circ
d_G) (x,y)\le -\kappa h_0\circ d_G(x,y), \ \forall (x,y)\in
S^2\backslash \triangle, \nneq which is equivalent to (by Theorem
\ref{thm21}): for the cost function $c(x,y)=h_0\circ d_G(x,y)$
$$
T_c(P_t(x,\cdot), P_t(y,\cdot)) \le c(x,y) e^{-\kappa t},\ t>0,
(x,y)\in S^2.
$$
To this purpose we begin by introducing a hypothesis which allows to
compare $d_G(X_t,Y_t)$  with some reference process $Z_t$ valued in
$\nn$, where $(X_t,Y_t)$ is the Markov process generated by a
coupling generator $\LL^\pi$ of $\LL$ with the kernel of jumping
rates $J^\pi((x,y), (x',y'))$.

Note that $J^{\pi}((x,y), (x',y'))=0$ if $|d_G(x',y')-d_G(x,y)|\ge
3$ for any coupling kernel $J^{\pi}$ of jumping  rates. We introduce
the following comparison condition on $J^\pi$:

\medskip

${\bf C(J;\alpha, \beta)}$: There exist two functions
$\alpha(x,y)\ge 1$ and
$\beta=(\beta_{-2},\beta_{-1},\beta_1,\beta_2): S^2\backslash
\triangle\to (\rr^+)^4$, and a nonnegative function $J_n(n+j)\ge 0$
on $-2\le j\le 2$ and $n\in [1,D_G]\cap\nn$ so that $0\le n+j\le
D_G$ (where $D_G=Diam(S,d_G)$ is the graph diameter of $S$, maybe
infinite), such that

\benu \item for every $x,y\in S$ with $d_G(x,y)=n\ge 1$, there is a
coupling jumping  rates measure $J^\pi((x,y), (x',y'))$ such that
for $j=1,2$,
 \beq\label{comparison1} \aligned
a_{\pi,j}(x,y)&:=\sum_{(x',y'):\ d_G(x',y')=n-j}J^\pi((x,y),
(x',y'))\\
&\ge \alpha(x,y)[
J_n(n-j)+\beta_{-j}(x,y)];\\
 b_{\pi,j}(x,y)&:=\sum_{(x',y'):\
d_G(x',y')=n+j}J^\pi((x,y), (x',y'))\\
&\le \alpha(x,y)[J_n(n+j)+\beta_j(x,y)];\\
\endaligned
 \nneq
 \item $(\beta_{-1}+2\beta_{-2})-(\beta_1+2\beta_2)\ge 0$ on $S^2\backslash \triangle$;
 \item $J_n(n-1)+2J_n(n-2)>0$ for all $1\le n\le D_G$. \nenu

Consider the generalized birth-death process $(Z_t)_{t\ge0}$ killed
at $0$ with generator \beq\label{comparison-model} \LL_{\rm ref}
f(n) =\sum_{j=-2}^2 J_n(n+j)[f(n+j)-f(n)], \ n\in [1,D_G]\cap\nn.
\nneq (This makes sense for $J_n(n+j)=0$ once if $n+j\notin
[0,D_G]\cap \nn$.) Letting $\tau_0:=\inf\{t\ge0; \ Z_t=0\}$ be the
first hitting time to $0$ of $Z_t$, then $Z_t=0$ for all $t\ge
\tau_0$, a.s. (as it is killed at $0$).

\bthm\label{thm-comparison} Assume ${\bf C(J;\alpha,\beta)}$.
Suppose that there is some increasing function
$h_0:[0,D_G]\cap\nn\to \rr$ with $h_0(0)=0$ and a positive constant
$\kappa>0$ such that $-\LL_{\rm ref} h_0(n)\ge \kappa h_0(n)$ for
all $n\in [1,D_G]\cap \nn$.

\benu[(a)]
\item
If $D_+h_0$ is non-increasing, then for the metric
$d(x,y):=h_0(d_G(x,y))$, $$Ric(\LL, d)\ge \kappa.$$

\item If $\beta(x,y)=0$ in  ${\bf C(J;\alpha,\beta)}$,
then
$$
\LL^\pi (h_0\circ d_G) (x,y)\le -\kappa (h_0\circ d_G) (x,y), \
(x,y)\in S^2,
$$
where  $\LL^\pi$ is the Markov generator associated with the jumping
rates kernel $J^\pi$ in ${\bf C(J;\alpha,\beta)}$. In particular for
the cost-function $c(x,y)=h_0(d_G(x,y))$,
$$
T_c(P_t(x,\cdot), P_t(y,\cdot))\le e^{-\kappa t} c(x,y),\ t\ge0,
(x,y)\in S^2.
$$
\nenu
 \nthm

\bprf (a) Since $h_0(0)=0$ and $D_+h_0>0$ is non-increasing,
$d(x,y)=h_0(d_G(x,y))$ is a metric.  If $d_G(x,y)=n\ge1$, we have by
${\bf C(J,\alpha, \beta)}$
$$\aligned
\LL^\pi (h_0\circ d_G)(x,y)&=\sum_{(x',y')} J^\pi((x,y),
(x',y'))[h_0(d_G(x',y')) - h_0(d_G(x,y))]\\
&\le \sum_{j=-2}^2 \alpha(x,y)(J_n(n+j)+\beta_j(x,y))[h_0(n+j)-h_0(n)]\\
&\le \alpha(x,y) \LL_{\rm ref} h_0(n)\le -\alpha(x,y)\kappa h_0(n)\\
&\le - \kappa h_0\circ d_G(x,y),
\endaligned$$
where the inequality in the third-line above  holds because
$$\sum_{j=-2}^2\beta_{j} [h_0(n+j)-h_0(n)]\le (\beta_1 +
2\beta_2) D_+h_0(n) -(\beta_{-1} + 2\beta_{-2}) D_+h_0(n-1)\le0
$$ by our condition (2) in ${\bf C(J,\alpha, \beta)}$ and $D_+h_0(n-1)\ge D_+h_0(n)$. Hence $Ric(\LL, d)\ge \kappa$ by Theorem
\ref{thm21}.

{\bf (b).} Now since $\beta(x,y)=0$, the argument above works
without the non-increasingness of $D_+h_0$ and gives us
$$
\LL^\pi h_0\circ d_G(x,y)\le  -\kappa \alpha(x,y) h_0\circ
d_G(x,y)\le  -\kappa  h_0\circ d_G(x,y).
$$
The conclusion in this part follows by Theorem \ref{thm21} again.
 \nprf

\brmk{\rm The best choice of $\kappa$ in the comparison theorem
above must be the smallest eigenvalue $\lambda_0$ of $-\LL_{\rm
ref}$ with the Dirichlet boundary condition at $0$ and $h_0$ the
associated positive eigenfunction (if exists). } \nrmk

\subsection{Comparison with a birth-death process} If $J_n(n+j)=0$
for $j=\pm 2$ in Theorem \ref{thm-comparison}, the reference process
generated by $\LL_{\rm ref}$ becomes a usual birth-death process for
which many problems admit explicit solutions (\cite{Chen92}). In
this paragraph we will provide some explicit estimates.

The following corollary yields an explicit quantitative estimate for
$$\ee_{(x,y)}\int_0^{\tau_c}  g(d_G(X_t,Y_t))dt,$$ where
$g:[0,D_G]\cap\nn\to \rr$ is nonnegative with $g(0)=0$, $(X_t,Y_t)$
is the Markov process starting from $(x,y)$, generated by some
coupling generator $\LL^\pi$ satisfying ${\bf C(J,\alpha, \beta)}$,
$\tau_c=\inf\{t\ge0; X_t=Y_t\}$ is the coupling time. Let
$$\mu_{\rm ref}(k)=\frac{J_1(2)J_2(3)\cdots
J_{k-1}(k)}{J_2(1)\cdots J_k(k-1)}, \ 1\le k\le D_G$$ be the
symmetric measure of $\LL_{\rm ref}$.

 \bcor\label{cor31} Assume that
${\bf C(J,\alpha, \beta)}$ is satisfied for $\beta=0$ by some
coupling Markov generator $\LL^\pi$ so that $J_n(n+j)=0$ for $j=\pm
2$. Then for every fixed positive function $g: [1,D_G]\to
(0,+\infty)$ such that $\sum_{k\ge 1} g(k)\mu_{\rm ref}(k)<+\infty$,
and any function $f$ on $S$ such that $|f(x)-f(y)|\le g(d_G(x,y))$
for all $(x,y)\in S^2$,  \beq\label{cor31a}\aligned \int_0^\infty
|P_tf(x)-P_tf(y)|dt
&\le \ee_{(x,y)}\int_0^{\tau_c} g(d_G(X_t,Y_t))dt\\
&\le \sum_{n=1}^{d_G(x,y)}\frac{\sum_{k\ge n} g(k)\mu_{\rm
ref}(k)}{\mu_{\rm ref}(n) J_n(n-1)}, \ \forall x\ne y(\in S)
\endaligned\nneq
where $(X_t,Y_t)$
is the Markov process starting from $(x,y)$, generated by some
coupling generator $\LL^\pi$ satisfying the comparison condition ${\bf C(J,\alpha, \beta)}$.

\ncor

\bprf  For every $(x,y)\in S^2\backslash \triangle$,
$$\aligned
\int_0^{+\infty} |P_tf(x)-P_tf(y)| dt&=\int_0^{+\infty}
|\ee_{(x,y)}[
f(X_t)-f(Y_t)]| dt\\
&\le \ee_{(x,y)} \int_0^{\tau_c} g(d_G(X_t,Y_t)) dt.
\endaligned
$$
Let $h:[0,D_G]\cap\nn\to \rr$ be the increasing function determined
by $h(0)=0$ and
$$
h(n)-h(n-1)=\frac{\sum_{k\ge n} g(k)\mu_{\rm ref}(k)}{\mu_{\rm
ref}(n) J_n(n-1)}, \ n\ge1.
$$
It is a solution to the Poisson equation $-\LL_{\rm ref} h(k)=g(k)$
for $k\ge 1$ (with the Dirichlet boundary condition $h(0)=0$). By
Theorem \ref{thm-comparison}(b) and its proof, if $n=d_G(x,y)\ge1$
$$
\LL^\pi (h\circ d_G)(x,y)\le  \alpha(x,y) (\LL_{\rm ref} h)(n)= -
 \alpha(x,y) g(n)\le -g\circ d_G (x,y).
$$
Setting $g(0):=0$, this inequality holds automatically if $x=y$ by
the definition of the coupling generator $\LL^\pi$.  Then
$$
(h\circ d_G)(X_t,Y_t) - h\circ d_G(x,y)+ \int_0^t g(d_G(X_s,Y_s))ds
$$
is a supermartingale. Therefore by Fatou's lemma,
$$\ee_{(x,y)}
\int_0^{\tau_c} g(d_G(X_t,Y_t)) dt\le h(d_G(x,y)),
$$
which yields the claim (\ref{cor31a}).
 \nprf

\brmk{\rm The estimate (\ref{cor31a}) for the solution $F=\int_0^\infty P_tf dt$ of the Poisson equation $-\LL F=f$ when $\mu(f)=0$ can be used
to obtain the transportation-information inequalities in A. Guillin {\it et al.} \cite{GLWY, GLWW,GJLW},
which are equivalent to the concentration inequalities of empirical means.

See Joulin and Ollivier \cite{JO10} for concentration inequalities for discrete time Markov chains under the positive Ricci curvature condition.
}\nrmk

We turn to give some explicit estimates in some typical cases.

\bcor\label{cor32} Assume that  ${\bf C(J,\alpha, \beta)}$  is
satisfied by some coupling Markov generator $\LL^\pi$ so that
$J_n(n+j)=0$ for $j=\pm 2$, and
$$
J_n(n-1)= a,\ J_n(n+1)= b
$$
for all $n\ge1$, for some positive constants $a,b$.

\benu[(a)]
\item If $a>b$ and $\beta(x,y)=0$ in  ${\bf C(J,\alpha, \beta)}$, then for $h_0(n)=1_{n>0} e^{cn}$ with $c=\frac 12 (\log a-\log
b)$ and $c(x,y)=h_0(d_G(x,y))$,
$$
T_c(P_t(x,\cdot), P_t(y,\cdot)) \le \exp\left(-(\sqrt{a} -
\sqrt{b})^2 t \right) h_0(d_G(x,y)).
$$
\item If $a=b$ and the graph diameter $D_G$ of $S$ is finite, then for the metric $d(x,y)=h_0(d_G(x,y))$ where $h_0(k)=\sin
\frac{k\pi}{2D_G}$,
$$
Ric(\LL, d)\ge 2a[1-\cos \frac{\pi}{2D_G}].
$$
\item If $a<b$ and $D_G<+\infty$, letting $h_0: [0,D_G]\cap \nn\to \rr$  be the increasing function
 determined by $h_0(0)=0$,
$$
h_0(n)-h_0(n-1)=\frac 1a \sum_{k=n}^{D_G} \left(\frac
ba\right)^{k-n} = \frac{\left(\frac ba\right)^{D_G-n+1}-1}{b-a},
$$
we have for the metric $d(x,y)=h_0(d_G(x,y))$,
$$
Ric(\LL, d)\ge \frac 1{h_0(D_G)} =
\frac{(b-a)^2}{b[(b/a)^{D_G}-1]-D_G(b-a)}.
$$

 \nenu
\ncor

\bprf (a) For the function $h_0$ given here, we have $\LL_{\rm ref}
h_0(n)=-(\sqrt{a}-\sqrt{b})^2 h_0(n)$ for all $1\le n<D_G$ and
$\LL_{\rm ref} h_0(D_G)< -(\sqrt{a}-\sqrt{b})^2 h_0(D_G)$ if
$D_G<+\infty$. It remains to apply Theorem \ref{thm-comparison}.

(b) For the function $h_0$ given here, we have $\LL_{\rm ref}
h_0(n)=-2a(1-\cos\frac{\pi}{2D_G}) h_0(n)$ for all $1\le n\le D_G$.
Since $h_0$ is concave on $k\in [0,D_G]$, $h_0\circ d_G$ is a
metric. Applying Theorem \ref{thm-comparison}, we conclude (b).

(c) For the function $h_0$ given in this part, we have $-\LL_{\rm
ref} h_0(n)=1$ for all $1\le n\le D_G$. Then
$$
-\LL_{\rm ref} h_0(n)\ge \frac 1{h_0(D_G)} h_0(n), \ 1\le n\le D_G.
$$
As $D_+h_0$ is decreasing, $d=h_0\circ d_G$ is again a metric. We
get (c) by Theorem \ref{thm-comparison}.
 \nprf

\brmk\label{cor32-rem}{\rm If the coupling generator $\LL^\pi$ is of the one-step type fixed
in the corollary above (i.e. $J_n(n+j)=0$ for $j=\pm 2$), the optimal estimate of the Ricci curvature lower bound $\kappa$
 should be the smallest eigenvalue $\lambda_0$ of $-\LL_{\rm ref}$. In that sense
  the estimate in part (a) is asymptotically optimal when
$D_G\to +\infty$, and the estimate in (b) is optimal because $2a[1-\cos\frac \pi{2D_G}]$ is
the smallest eigenvalue $\lambda_0$ of $-\LL_{\rm ref}$.
However the estimate of the Ricci curvature lower bound $\kappa$
in part (c) is not sharp: the optimal one
 should be the smallest eigenvalue $\lambda_0$ of $-\LL_{\rm ref}$. Let $h$
 be the positive eigenfunction associated with $\lambda_0$ (Perron-Frobenius theorem), with
 $h(0)=0$. Then
 $$
h(n)-h(n-1)=\frac {\lambda_0}a \sum_{k=n}^{D_G} (b/a)^{k-n} h(k),
 $$
which implies that $h$ is increasing. Then
$$
h(1) \ge \frac {\lambda_0}{a} \sum_{k=1}^{D_G} (b/a)^{k-1} h(1),
$$
which yields
$$
\lambda_0\le \frac{b-a}{(b/a)^{D_G} - 1 }.
$$
This shows that the estimate of $\kappa$ in part (c) is of correct
order $(b/a)^{-D_G}$ in $D_G$. }\nrmk

\subsection{Estimate of eigenvalues in terms of degree and diameter}
The following result is an improvement of Lin-Yau \cite[Theorem
1.8]{LY10} and generalizes their result to non-symmetric case.

 \bcor\label{cor33}
For the generator $\LL$ on a finite graph $(S,E)$ with the graph
diameter $D_G\ge 2$,  for any eigenvalue $\lambda\ne 0$ of $-\LL$ in
$\cc$,
\beq\label{cor33a}\aligned {\rm Re}(\lambda)&\ge
\frac{2\lambda_*(d_\LL-2)}{d_\LL\left[\sum_{k=1}^{D_G}(d_\LL-1)^k-D_G\right]},\endaligned\nneq
where $\lambda_*:=\inf_{x\in S}\lambda(x)$ is minimal
jumping rate (recalling that $\lambda(x)=\sum_{y\in S: y\sim
x}J(x,y)$ is the total jumping rate at $x$) and \beq\label{cor33b}
d_\LL=\sup_{x\in S}\max_{y\sim x} \frac{\lambda(x)}{J(x,y)}\nneq is
the maximal degree w.r.t. $\LL$, assumed to be $>2$. \ncor

As $d_\LL\ge \sup_{x\in S} d_x$ (by taking $y\sim x$ so that
$J(x,y)=\min_{x'\sim x}J(x,x')$ in (\ref{cor33b})), then once if
$d_\LL\le 2$, $(S,E)$ can be identified as a finite interval of
$\nn$ (for which we have Chen's variational formula of $\lambda_1$) or
as a discrete circle (which will be studied below).

\bprf Notice that in this finite space case, for any eigenvalue
$\lambda\ne 0$ of $-\LL$, $ {\rm Re}(\lambda)\ge \inf_{(x,y)\in
S^2\backslash \triangle} Ric_{(x,y)}(\LL, d)$ for any metric $d$. Let us
find such a metric $d$ by means of Theorem \ref{thm-comparison}.

We use the independent coupling generator $\LL^\pi$ with jumping
rates kernel $J^\pi$ given in (\ref{indepC}). Remark that for any
$y\sim x$,
$$
J(x,y)\ge \frac{\lambda(x)}{d_\LL}.
$$
For every $(x,y)\in S^2$ with $d_G(x,y)=n\ge1$, there is at least
one geodesic from $x$ to $y$. If $x_1(\sim x),y_1(\sim y)$ are in
this geodesic, then $d(x_1,y)=d(y_1,x)=n-1$, and
$$\aligned
J^{\pi}((x,y),(x_1,y))+ J^{\pi}((x,y),(x,y_1))&=J(x,x_1)+J(y,y_1)\\
&\ge
\frac{\lambda(x)}{d_\LL} + \frac {\lambda(y)}{d_\LL}\ge \frac{2\lambda_*}{d_\LL}\cdot \frac{\lambda(x)+\lambda(y)}{2\lambda_*},\\
\sum_{(x',y'): d_G(x',y')=n+1}J^{\pi}((x,y),(x', y')) &\le
\lambda(x)
-J(x,x_1)+\lambda(y)-J(y,y_1)\\
&\le \lambda(x) - \frac{\lambda(x)}{d_\LL} +\lambda(y)-
\frac{\lambda(y)}{d_\LL} \\
&\le \frac{2\lambda_*(d_\LL-1)}{d_\LL}\cdot
\frac{\lambda(x)+\lambda(y)}{2\lambda_*}.
\endaligned$$ Therefore
the condition {\bf ${\bf C(J,\alpha, \beta)}$} is satisfied for
$\alpha(x,y)=\frac{\lambda(x)+\lambda(y)}{2\lambda_*}\ge1,\
\beta(x,y)=0$ and
$$J_n(n-1)=\frac {2\lambda_*}{d_\LL} (1\le n\le D_G), \ J_n(n+1)=\frac{2\lambda_*(d_\LL-1)}{d_\LL}\ (0\le n\le D_G-1).$$
Let $h_0:[0,D_G]\cap \nn\to\rr^+$ be a positive eigenfunction of
$-\LL_{\rm ref}$ with the Dirichlet boundary condition at $0$ (so
$h_0(0)=0$), associated with the smallest eigenvalue $\lambda_0(\LL_{\rm
ref})$ (Perron-Frobenius theorem), which is increasing by Remark \ref{cor32-rem}. For the metric $d(x,y)=h_0(d_G(x,y))$, by Theorem
\ref{thm-comparison}, $Ric(\LL, d)\ge \lambda_0(\LL_{\rm ref})$. But
applying Corollary \ref{cor32}(c) with $a=2\lambda_*/d_{\LL}<b=2\lambda_*(d_\LL
-1)/d_\LL$, $\lambda_0(\LL_{\rm ref})$ is bounded from below by the
r.h.s. of (\ref{cor33a}).
 \nprf

\brmk{\rm For the special case of the Laplace operator on a finite
graph, a beautiful lower bound only in terms of the  degrees
 is due to Barlow, Coulhon and Grigoryan
\cite{BCG01} and Chung \cite{Chung97}:
$$
\lambda_1(\Delta) \ge \frac 1{d_* |E|},
$$
where $d_*=\max_{x\in S}d_x$, $|E|=\sum_{x\in S} d_x$ is the number
of oriented edges in $E$. In \cite[Section 6.5]{Chung97},
given $d\ge3$ and $D_G\ge 2$, a $d$-regular graph of diameter $D_G$ can
have as many as $d(d-1)^{D_G}$ vertices, then $|E|$ is
$d^2(d-1)^{D_G}$. The lower bound above becomes $\lambda_1(\Delta)
\ge \frac 1{d^3(d-1)^{D_G}}$, which is of the same order
$(d-1)^{-D_G}$ as in the corollary above for big $D_G$.

Lin and Yau \cite[Theorem 1.8]{LY10} generalized the above result of
\cite{BCG01, Chung97} to the case where $\lambda(x)\equiv 1$,
$\mu(x)J(x,y)=\mu(y)J(y,x)$ (the symmetry of $\LL$ on $L^2(\mu)$;
this quantity is $\mu_{xy}$ in \cite{LY10}) and showed that the
spectral gap of $\LL$ satisfies
\beq
\lambda_1\ge \frac{1}{d_\LL D_G \exp\left(d_\LL D_G+1\right)-1}.
\nneq
Our lower bound (\ref{cor33a}) is better.
 } \nrmk

 \brmk{\rm For estimates of the spectral gap $\lambda_1$ on graphes by means of other tools,
 the reader is referred to Lawler and Sinclair \cite{LS88}, Diaconis and Stroock \cite{DiSt91}, and to the book
 of Chung \cite{Chung97} for references.  See  Ma {\it et al.} \cite{MWW15}, Liu {\it et al.} \cite{LMW}  for some recent progresses.
 }\nrmk

\subsection{One step optimal coupling}

The main difference of graphs from Riemmanian manifolds is: the
geodesics in the graph metric $d_G$ linking two vertices are in
general not unique. Given a Markov generator $\LL$,  there are often
many optimal coupling generators of $\LL$ in the metric $d_G$,
because of non-uniqueness of geodesics.

In the following result we
show that for the Ricci curvature $Ric_{(x_0,y_0)}:=Ric_{(x_0,y_0)}(\LL, d_G)$ w.r.t. the graph metric, we can always choose an one-step optimal coupling $\LL^\pi$  so that our previous comparison results are
applicable.

\blem\label{lem33} There is always a $d_G$-optimal coupling
$\LL^\pi$ of $\LL$, i.e. \beq\label{lem33c} \LL^\pi d_G(x_0,y_0)=
-Ric_{(x_0,y_0)}\cdot d_G(x_0, y_0), \ (x_0,y_0)\in S^2\backslash
\triangle, \nneq where $Ric_{(x_0,y_0)}:=Ric_{(x_0,y_0)}(\LL, d_G)$,
so that its coupling kernel of jumping rates $J^\pi$ satisfies: for
every $(x_0,y_0)\in S^2\backslash \triangle$,

\benu
\item
\beq\label{lem33a} \sum_{(x',y'): d_G(x',y')=d_G(x_0,y_0)\pm 2}
J^\pi((x_0,y_0), (x',y'))=0; \nneq
\item for any neighbor $x_1$ of $x_0$ and neighbor $y_1$ of $y_0$
both lying to a geodesic linking $x_0$ to $y_0$ ($x_1,y_1$ may
be the same vertex if $d_G(x_0,y_0)=2$), \beq\label{lem33b}
J^\pi((x_0,y_0), (x_1,y_0))=J(x_0,x_1),\ J^\pi((x_0,y_0),
(x_0,y_1))=J(y_0,y_1).\nneq \nenu
 \nlem

The specification (\ref{lem33b}) means that when one goes from $x_0$ to $x_1$ (closer to $y_0$) in the direction of the geodesic, the other stays at $y_0$, and {\it vice versa}. This choice is not the good one for getting together more rapidly, but well adapted for comparison.

\bprf Let $J^\pi((x_0,y_0),\cdot)$ be an optimal coupling in the
Wasserstein transport cost $T_{d_G}$ of $J(x_0,\cdot), J(y_0,\cdot)$
with $J(x_0,x_0)=\lambda(y_0)$ and $J(y_0,y_0)=\lambda(x_0)$,
$\LL^\pi$ the corresponding coupling generator of $\LL$. By Theorem
\ref{thm22} (and in the re-definition \ref{def_Ricci}), $\LL^\pi$
satisfies (\ref{lem33c}).

{\bf Step 1. }  At first we show that we can construct a new optimal
coupling generator $\LL^{\pi'}$, so that (\ref{lem33a}) is
satisfied.  Define
$$\aligned
&\LL^{\pi'}F(x_0,y_0)=\sum_{(x',y'): |d_G(x',y')-d_G(x_0,y_0)|\le
1}
[F(x',y')-F(x_0,y_0)]J^{\pi}((x_0,y_0), (x',y'))\\
&+\sum_{(x',y'): d_G(x',y')=d_G(x_0,y_0)\pm 2} [F(x',y_0)+
F(x_0,y')-2F(x_0,y_0)]J^{\pi}((x_0,y_0), (x',y')).
\endaligned
$$
This is again a coupling generator of $\LL$, which replaces the jump
from $(x_0,y_0)$ to $(x',y')$ so that $d_G(x',y')=d_G(x_0,y_0)\pm 2$
with rate $J^{\pi}((x_0,y_0), (x',y'))$ by two free jumps from
$(x_0,y_0)$ to $(x',y_0)$ or to $(x_0,y')$, both with rate
$J^{\pi}((x_0,y_0), (x',y'))$.

We see that
$$\aligned
\sum_{(x',y'): d_G(x',y')=d_G(x_0,y_0)+ 1} & J^{\pi'}((x_0,y_0),
(x',y'))=\sum_{(x',y'): d_G(x',y')=d_G(x_0,y_0)+ 1}
J^{\pi}((x_0,y_0),
(x',y'))\\
&+ 2 \sum_{(x',y'): d_G(x',y')=d_G(x_0,y_0)+ 2} J^{\pi}((x_0,y_0),
(x',y'))
\endaligned$$
and
$$\aligned
\sum_{(x',y'): d_G(x',y')=d_G(x_0,y_0)- 1} &J^{\pi'}((x_0,y_0),
(x',y'))=\sum_{(x',y'): d_G(x',y')=d_G(x_0,y_0)-1}
J^{\pi}((x_0,y_0),
(x',y'))\\
&+ 2 \sum_{(x',y'): d_G(x',y')=d_G(x_0,y_0)- 2} J^{\pi}((x_0,y_0),
(x',y')).
\endaligned$$
Therefore $\LL^{\pi'} d_G(x_0,y_0)= \LL^\pi d_G(x_0,y_0)$:
$\LL^{\pi'}$ is again an optimal coupling of $\LL$ for the graph
metric $d_G$.

{\bf Step 2.} By Step 1 we may and will assume that
$J^\pi((x_0,y_0),\cdot)$ satisfies (\ref{lem33a}).

Let $d_G(x_0,y_0)=n\ge1$ and  consider two vertices  $x_1\sim x,
y_1\sim y$ in a geodesic $\gamma_{x_0,y_0}$ linking $x_0$ to $y_0$.
By the definition (\ref{21}) of coupling and (\ref{lem33a}),
\beq\label{lem33d}\begin{cases} J(x_0,x_1)&=J^\pi((x_0,y_0), (x_1,
y_0))+\sum_{y'\sim y_0: y'\ne
y_1} J^\pi((x_0,y_0), (x_1, y'));\\
J(y_0,y_1)&=J^\pi((x_0,y_0), (x_0, y_1))+\sum_{x'\sim x_0: x'\ne
x_1} J^\pi((x_0,y_0), (x', y_1)) . \end{cases}\nneq  Define the new
coupling
$$
J^{\pi'}((x_0,y_0), (x',y')):= \begin{cases} J(x_0,x_1),\ &\text{ if
}\ (x',y')=(x_1,y_0);\\
0, &\text{ if }\ x'=x_1,\ y'\sim y_0;\\
J(y_0,y_1),\ &\text{ if
}\ (x',y')=(x_0,y_1);\\
0, &\text{ if }\ x'\sim x_0,\ y'=y_1;\\
J^\pi((x_0,y_0), (x',y')), &\text{ otherwise. }
\end{cases}
$$
$J^{\pi'}((x_0,y_0),\cdot)$ satisfies again (\ref{lem33d}): it is
again a coupling of $J(x_0,\cdot), J(y_0,\cdot)$. For any neighbor
$y'$ of $y_0$ different from $y_1$ such that $J^\pi((x_0,y_0), (x_1,
y'))>0$, $d_G(y',x_1)=n$ or $n-1$ (the possibility of $n-2$ is
excluded by (\ref{lem33a})); and for any neighbor $x'$ of $x_0$
different from $x_1$ such that $J^\pi((x_0,y_0), (x', y_1))>0$,
$d_G(x',y_1)=n$ or $n-1$. Thus
$$\sum_{(x',y')\in S^2} d_G(x', y') J^{\pi'}((x_0,y_0), (x',y'))\le
\sum_{(x',y')\in S^2} d_G(x', y') J^{\pi}((x_0,y_0), (x',y'))$$ and thus
$\LL^{\pi'}d_G(x_0,y_0)\le \LL^\pi d_G(x_0,y_0)$. As $\LL^\pi$ is an
optimal coupling, so is $\LL^{\pi'}$. And $J^{\pi'}$ verifies
(\ref{lem33b}).
 \nprf

\subsection{Discrete cycle}
Now we present several examples.
We begin with the over-worked model: random walk on the discrete
circle.

\bexa\label{discrete cycle example}{\bf (discrete circle) }{\rm Consider the Laplacian $\Delta$ on
 $S=\zz/n\zz$ $(n\ge3)$ which can be identified
as the discrete circle $\{s_k=\exp(i\frac{2k\pi}{n}); 0\le k\le
n-1\}$. Its graph diameter is $D_G=[n/2]$, the integer part of
$n/2$. By Fourier Analysis, the spectrum of $-\Delta$ is
$\{\lambda_k=1-\cos \frac{2k\pi}{n};\ 0\le k\le n-1\}$. The
logarithmic Sobolev inequality with sharp constant $\lambda_1$ was
established by Chen and Sheu \cite{CS03}.
See Sammer and Tetali \cite{ST09} for concentration inequalities on the torus (the product space of discrete cycles).

If one uses the independent coupling, the condition of  Corollary
\ref{cor32} is satisfied with $a=b=1$. Then for $h_0(k)=\sin
\frac{k\pi}{2D_G}$ and $d(x,y)=h_0(d_G(x,y))$, $Ric(\Delta, d)\ge
2(1-\cos \frac{\pi}{2D_G})$. This lower bound is asymptotically
equivalent to $\lambda_1/2$ for big $n$.

Now we will use a mixing of the reflection coupling and of the
independent one when $x,y$ are neighbors and the reflection coupling
otherwise to get the sharp result below.

 \bprop\label{prop-circle} For the metric
$$
d(x,y)=\sin \frac{d_G(x,y) \pi}{n}
$$
on the discrete circle $S=\zz/n\zz$ with $n\ge3$, we have
$Ric(\Delta, d)\ge \lambda_1=1-\cos \frac{2\pi}{n}$. \nprop

\bprf {\bf Step 1.} We begin by the construction of the coupling. If
$x,y$ are neighbors, say $y=x+1$,  the coupling jumping rates kernel
$J^\pi((x,y),\cdot)$ will be given by
$$
J^\pi((x,y), (x',y'))=\begin{cases}
\frac 12, \ &\text{ if  } (x',y')=(y,y);\\
\frac 12, &\text{ if  } (x',y')=(x,x);\\
\frac 12, &\text{ if   } (x',y')=(x-1,y+1);\\
0, &\text{ otherwise }.
\end{cases}
$$
One can think this coupling in the following way: when the walker A
at position $x$ goes to $y=x+1$ at the next step, he does not inform
his co-walker B at position $y$ and {\it vice-versa} (independent
coupling); but if he goes to $x-1$, he informs his friend who will
go simultaneously at the opposite direction $y+1$ (but if $n=3$,
$x-1=y+1(\mbox{mod}\ 3)$, contrary to their attentions): that is the
reflection coupling.

Now if $d_G(x,y)\ge 2$ (necessarily $n\ge 4$), we take the
reflection coupling given by
$$
\LL^\pi f(x,y)=\ee f(x+\xi, y-\xi) - f(x,y),
$$
where $\xi$ is a random variable with $\pp(\xi=1)=\pp(\xi=-1)=1/2$.
For this coupling,  we see that $d_G(X_t,Y_t)$ is a Markov process
whose generator coincides with $\LL_{\rm ref}$ for $J_k(k+j)$ given
below:

\benu \item  for $k=1$,
$$
J_1(1+j)=\begin{cases}
1, \ &\text{ if } j=-1;\\
\frac 12, &\text{ if   } j= 2;\\
0, &\text{ otherwise;}
\end{cases}
$$
\item for $2\le k\le D_G-2$, $J_k(k+j)=\frac 12$ for $j=\pm 2$ and
$J_k(k+j)=0$ for $j=\pm 1$;

\item finally for $k=D_G-1$ or $D_G$:  if $n$ is even (then $n=2D_G$),
$$J_k(k+j)=\begin{cases}
\frac 12, &\text{ if }  k=D_G-1, j=-2 \mbox{ or } j=0;\\
1, &\text{ if  } k=D_G,j=-2;\\
0, &\text{ otherwise}
\end{cases}
$$
and if $n$ is odd, i.e. $n=2D_G+1$,
$$J_k(k+j)=\begin{cases}
\frac 12, &\text{ if }  k=D_G-1, j=-2 \mbox{ or } j=1;\\
\frac 12, &\text{ if  } k=D_G,j=-1 \mbox{ or } -2;\\
0, &\text{ otherwise}.
\end{cases}
$$
\nenu Therefore the comparison condition ${\bf C(J,\alpha,\beta)}$
is satisfied with $J_k(k+j)$ given above, $\alpha(x,y)=1$,
$\beta(x,y)=0$.

{\bf Step 2.} Notice that $d(x,y)=h(d_G(x,y))$ ($x,y\in S$) where
$$
h(x)=\sin \frac{x\pi}{n}, \ x\in\rr
$$
satisfies: for all $x\in\rr$, \beq\label{circle2}
\frac{1}{2}\{[h(x-2)-h(x)]+[h(x+2)-h(x)]\} = -\lambda_1 h(x), \
h(x)=h(n-x). \nneq By Step 1,
$$
(\LL^\pi d)(x,y)=(\LL^\pi h\circ d_G ) (x,y)=(\LL_{\rm ref}
h)(d_G(x,y)),
$$
where $\LL_{\rm ref} h(k)=\sum_{j=-2}^2[h(k+j)-h(k)] J_k(k+j)$ with
jumping rates $J_k(k+j)$ given in Step 1. For the conclusion of this
proposition, by Theorem \ref{thm-comparison}, it suffices to prove
\beq\label{circle3} \LL_{\rm ref} h(k)=-\lambda_1 h(k), \ 1\le k\le
D_G. \nneq

At first the equality above holds for $2\le k\le D_G-2$,  by
(\ref{circle2}).

For $k=1$, noting that $h(-1)=-h(1)$ and $h(0)=0$, we have by
(\ref{circle2}),
$$\aligned
\LL_{\rm ref} h(1)&=\frac 12 [h(3)-h(1)] + [h(0)-h(1)]\\
&= \frac 12 \{[h(3)-h(1)] + [h(-1)-h(1)]\}=-\lambda_1 h(1).
\endaligned $$ For $k=D_G-1$ or $D_G$, we separate our discussion into
two cases: $n$ is even or odd.

{\bf Case 1. $n$ is even. } For $k=D_G-1$, as $h(D_G+x)=h(D_G-x)$
(by (\ref{circle2})), we have
$$\aligned
\LL_{\rm ref} h(D_G-1)&=\frac 12[h(D_G-3)-h(D_G-1)]\\
&=\frac 12\{[h(D_G-3)-h(D_G-1)]+ [h(D_G+1)-h(D_G-1)]\}\\
&=-\lambda_1 h(D_G-1). \endaligned$$ For $k=D_G$, as
$h(D_G+2)=h(D_G-2)$,
$$\aligned
\LL_{\rm ref} h(D_G)&=[h(D_G-2)-h(D_G)]\\
&=\frac 12\{[h(D_G-2)-h(D_G)]+ [h(D_G+2)-h(D_G)]\}\\
&=-\lambda_1 h(D_G).  \endaligned$$

{\bf Case 2. $n$ is odd. } For $k=D_G-1$, by (\ref{circle2}) we have
$h(D_G)=h(D_G+1)$ and then
$$\aligned
\LL_{\rm ref} h(D_G-1)&=\frac 12\{[h(D_G-3)-h(D_G-1)]+[h(D_G)-h(D_G-1)\}\\
&=\frac 12\{[h(D_G-3)-h(D_G-1)]+[h(D_G+1)-h(D_G-1)]\}\\
&=-\lambda_1 h(D_G-1).\endaligned$$

Finally for $k=D_G$, as $h(D_G+2)=h(D_G-1)$,
$$\aligned
\LL_{\rm ref} h(D_G)&=\frac12\{[h(D_G-2)-h(D_G)]+[h(D_G-1)-h(D_G)]\}\\
&=\frac 12\{[h(D_G-2)-h(D_G)]+ [h(D_G+2)-h(D_G)]\}\\
&=-\lambda_1 h(D_G).  \endaligned$$ So we have completed  the proof
of (\ref{circle3}).  \nprf

 }\nexa

\subsection{Two-coloured graph}
\bexa {\bf (Two-partite complete graph or two-coloured graph)}
{\rm
Consider the Laplace generator $\LL=\Delta $ on $S$, where $S= S_1 \cup S_2$ is a two-partite complete graph i.e. a graph whose vertex set can be decomposeded into two disjoint and independent sets $S_1, S_2$ (no edge inside $S_j$),  with $|S_1|=N_1 \ge 2$, $|S_2|=N_2 \ge 2$, and every pair of vertices $x,y$ with $x \in S_1 $ and $y \in S_2$ is connected by an edge, in other words $E=(S_1\times S_2)\cup (S_2\times S_1)$. Its graph diameter is $D_G=2$.

 We construct the following coupling.

If $x\in S_1,y\in S_2$ or if $x\in S_2,y\in S_1$, we take the locally independent coupling jumping rates $ J^\pi((x,y), \cdot)$ at $(x,y)$.

If $x,y \in S_1 $ (resp. $S_2$), the coupling jumping rates kernel $ J^\pi((x,y), \cdot)$ will be given by
$$
J^\pi((x,y), (z,z))=\frac{1}{N_2}\  (\mbox{resp.}\  \frac{1}{N_1}), \quad  \mbox{if } z \in S_2 \  (\mbox{resp. } S_1), \ \ J^\pi((x,y), \triangle^c)=0.
$$

For this coupling, the reference generator $\LL_{\rm ref}$ for $J_k(k+j)$ is given below:
$$
\begin{aligned}
&J_1(0)=\frac{1}{N_1}+\frac{1}{N_2};\quad
J_1(2)=\frac{N_1-1}{N_1}+\frac{N_2-1}{N_2};\\
&J_2(0)=1;\quad
J_k(k+j)=0, \mbox{ otherwise}.
\end{aligned}
$$
We have $\LL^\pi h\circ d_G(x,y)=(\LL_{\rm ref} h)\circ d_G (x,y) $ for any function $h$ on $\{0,1,2\}$ with $h(0)=0$.  Setting
$$
h_0(0)=0, h_0(1)=1, \  0<h_0(2)\le \frac{N_1N_2}{2N_1N_2 -N_1-N_2}\le 1,
$$
we get by calculation $\LL_{\rm ref} h_0(k)\le - h_0(k), k=1,2$ and then
$$
\LL^\pi h_0\circ d_G(x,y)=(\LL_{\rm ref} h_0)\circ d_G (x,y) \le - h_0\circ d_G(x,y),\ x\ne y.
$$
Even if $h_0$ is not increasing, we check easily $d(x,y):=h_0\circ d_G(x,y)$ is again a metric. By Theorem \ref{thm21}, we have $Ric(\LL, d)\ge 1$.

For this model, letting $P=J$ be the transition probability kernel, we have $P^2(x,y)=\frac 1{N_j}$ if $x,y\in S_j$ and $P^2(x,y)=0$ otherwise. Hence the eigenvalues of $P^2$ are $0,1$. Therefore $\lambda_1=1$ (the spectral gap of $\Delta$), which shows that our esimate of the Ricci curvature lower bound above is optimal.

}\nexa

\subsection{Regular $k$-coloured graph}
\bexa {\bf (regular $k$-partite complete graph)}
{\rm Consider the regular $k(\ge 2)$-partite complete graph $(S,E)$, i.e. $S$ can be decomposed into $k$ disjoint parts $S_1,\cdots, S_k$ with $|S_i|=N\ge 2$ for any $i$, and
$(x,y)\in E$ if and only if $x\in S_i, y\in S_j$ for some $i\ne j$. Its graph diameter is $D_G=2$.

Consider the Laplace generator $\LL=\Delta $ on $S$.
We begin by the construction of the coupling.

If $x,y$ belong to some same coloured part $S_i $, the coupling jumping rates kernel $ J^\pi((x,y), \cdot)$ will be given by
$$
J^\pi((x,y), (z,z))=J(x,z)=J(y,z)=\frac{1}{N(k-1)}, \quad  \mbox{if } z \in S_j, j\neq i.
$$

If $x,y$ belong to two different coloured parts, i.e. $x\in S_i,y\in S_j, i\neq j$, the coupling jumping rates kernel $ J^\pi((x,y), \cdot)$ will be given by
$$
J^\pi((x,y), (x',y'))=
\begin{cases}
J(x,x')=\frac{1}{N(k-1)}, \quad & \mbox{if } x' \in S_j, y'=y; \\
J(y,y')=\frac{1}{N(k-1)}, \quad & \mbox{if } y' \in S_i, x'=x; \\
\frac{1}{N(k-1)}, \quad & \mbox{if }  x'=y'=z\in S_l,l\neq i, l\neq j;\\
0, & \mbox{ otherwise }.
\end{cases}
$$

For this coupling, the reference generator $\LL_{\rm ref}$ for $J_k(k+j)$ is given below:
$$
J_1(0)=\frac{2}{N(k-1)}+\frac{k-2}{k-1};\quad
J_1(2)=\frac{2(N-1)}{N(k-1)};\quad
J_2(0)=1;\quad
J_k(k+j)=0, \mbox{  otherwise}.
$$
And we have $\LL^\pi h\circ d_G(x,y)=(\LL_{\rm ref} h)\circ d_G (x,y) $ for any function $h$ on $\{0,1,2\}$ with $h(0)=0$.  Letting
$$
h_0(0)=0, h_0(1)=1, \  0<h_0(2)\le \frac{N}{2(N-1)}\le 1,
$$
we get by calculation $\LL_{\rm ref} h_0(k)\le - h_0(k), k=1,2$ and then
$$
\LL^\pi h_0\circ d_G(x,y)=(\LL_{\rm ref} h_0)\circ d_G (x,y) \le - h_0\circ d_G(x,y),\ x\ne y.
$$
$d(x,y):=h_0\circ d_G(x,y)$ is again a metric. By Theorem \ref{thm21}, we have $Ric(\LL, d)\ge 1$.
}
\nexa

\section{Exponential convergence in $W_{1,d_G}$}
\subsection{A general result on the exponential convergence in $W_{1,d_G}$}
A. Eberle \cite{Eberle16}, Luo and Wang \cite{LuoWang16} proved that for the diffusion $dX_t=\sqrt{2} dB_t +b(X_t) dt$ in $\rr^d$, if the drift $b(x)$ is dissipative at infinity
\beq\label{dissipativity3}
\left\<\frac {x-y}{|x-y|}, b(x)-b(y)\right\> \le -\delta |x-y| + C 1_{|x-y|\le R} |x-y|,\ x,y\in\rr^d,
\nneq
where $\delta>0, C\ge0$ are two constants, then $W_1(P_t(x,\cdot), P_t(y,\cdot))\le K e^{-\kappa t}|x-y|$ (for some constants $\kappa>0, K\ge 1$), i.e.
its transition semigroup $(P_t)$ converges exponentially rapidly to its unique invariant probability measure $\mu$.

The following is the counterpart of their result on graphs.

\bthm\label{thm33} Assume that $Ric_{(x,y)}(\LL, d_G)$ is bounded from
below and positive for big $d_G(x,y)$, i.e. there are constants
$N\in\nn^*$, $\kappa_\infty>0$ and $R\ge 0$ such that
\beq\label{thm33a} Ric_{(x,y)}(\LL, d_G) \ge
\begin{cases}\kappa_\infty,\ &\text{ if }\
d_G(x,y)\ge N;\\
-R, &\text{ if }\ d_G(x,y)<N.\end{cases} \nneq If moreover
\beq\label{thm33b} J_*:=\inf_{(x,y)\in E} J(x,y)>0, \nneq then there
are constants $K\ge 1, \delta>0$ explicitly computable, such that
\beq\label{thm33c} W_{1,d_G}(P_t(x,\cdot), P_t(y,\cdot))\le K
e^{-\delta t} d_G(x,y),\ t\ge0, x,y\in S. \nneq \nthm

Notice that if $\lambda_*:=\inf_{x\in}\lambda(x)>0$ and the maximal
degree $d_\LL<+\infty$, $J(x,y)\ge \frac{\lambda_*}{d_\LL}$, the
condition (\ref{thm33b}) is verified.

\bprf Since (\ref{thm33a}) still holds for bigger $N$, we may assume
without lose of generality that $\kappa_\infty\cdot N\ge 2J_*$. We write
$Ric_{(x,y)}:=Ric_{(x,y)}(\LL, d_G)$ for simplicity of notation.

Let $\LL^\pi$ be an optimal coupling generator of $\LL$ w.r.t.
$d_G$, with the coupling kernel of jumping rates
$J^\pi((x,y),\cdot)$ satisfying (\ref{lem33a}) and (\ref{lem33b}),
for all $(x,y)\in S^2\backslash \triangle$, constructed in Lemma \ref{lem33}.

For $(x,y)\in S^2$ with $d_G(x,y)=n\ge 1$, let
$$\aligned a_\pi(x,y)&:=\sum_{d_G(x',y')=n-1}
J^\pi((x,y),(x',y'));\\
b_\pi(x,y)&:=\sum_{d_G(x',y')=n+1} J^\pi((x,y),(x',y')),
\endaligned$$
then $\LL^\pi d_G(x,y) = b_\pi(x,y)-a_\pi(x,y)=-Ric_{(x,y)}\cdot n$.
By (\ref{lem33b}) and our condition (\ref{thm33b}),
$$a_{\pi}(x,y)\ge J(x,x_1) + J(y,y_1)\ge 2J_*,$$
where $x_1\sim x$ and $y_1\sim y$ belong to a geodesic linking $x$
to $y$. Thus by our condition on the Ricci curvature
$$\beta_{-1}(x,y):=a_\pi(x,y)-\max\{2J_*, \kappa_\infty d_G(x,y)
1_{d_G(x,y)\ge N} \}\ge 0.
 $$
If $d_G(x,y)=n\ge 1$, $ a_\pi(x,y)=\beta_{-1}(x,y)+J_{n}(n-1)$,
where
$$
J_n(n-1)=\begin{cases} 2J_*,\ &\text{ if }\ n\in [1,N-1];\\
 \kappa_\infty n, &\text{ if }\ n\ge N.
\end{cases}
$$
On the other hand
$$\aligned
b_\pi(x,y)&=a_{\pi}(x,y)-Ric_{(x,y)}\cdot n \\
&=\beta_{-1}(x,y) + \max\{2J_*, \kappa_\infty n 1_{n\ge N}\}
-Ric_{(x,y)}\cdot n\\
&\le \beta_{-1}(x,y) +J_n(n+1),
\endaligned$$
where
$$
J_n(n+1)=\begin{cases} 2J_* +R n,\ &\text{ if }\ n\in [1,N-1];\\
 0, &\text{ if }\ n\ge N.
\end{cases}
$$
In other words the comparison condition ${\bf C(J,\alpha, \beta)}$
is satisfied for $\alpha(x,y)=1$, $\beta_1(x,y)=\beta_{-1}(x,y)$ and
$J_n(n\pm 1)$ given above, and $J_n(n+j)=\beta_j=0$ for $j=\pm 2$.

Let
$$g(n)=\begin{cases}\kappa_\infty\cdot n, \ &\text{ if }\ n\ge N+1;\\
2J_* &\text{ if }\ n\in [1,N]. \end{cases}$$ The solution $h:\nn\to
\rr$ of the Poisson equation
$$
-\LL_{\rm ref} h(n)=g(n), n\in \nn^*
$$
so that $h(k)=k$ for $k\ge N$ (in fact the equation above is
verified for $n\ge N+1$), is determined by
$$
D_+h(n-1)= \frac{\nu[n,N]}{\nu(n)},\ n\in [1,N],
$$
where
\beq\label{thm33e}
 \nu(N)=1,\ \nu(n) =
\frac{(2J_*)^{N-n}}{\prod_{k=n}^{N-1} (2J_* +Rk)}, \ n\in [1,N-1] \text{ and }
\nu[n,N]=\sum_{k=n}^N \nu(k).
\nneq
As $\nu(n+1)=\frac{2J_*+ Rn}{2J_*} \nu(n)\ge \nu(n)$,
$D_+h(n)$ is decreasing in $n\in [1,N]$ and then over $\nn^*$ as
$D^+h(N-1)=1=D_+h(k)$ for all $k\ge N$.

Finally let
$$\aligned h_0(n)&:=h(n)-h(0)= \sum_{k=1}^n D_+ h(k-1)\\
&= \sum_{k=1}^n \sum_{j=k}^N  \frac{\nu[j,N]}{\nu(j)} +(n-N)^+,
\endaligned $$ which is increasing and $D_+h_0(n)$ is decreasing. We have
$$
\LL_{\rm ref} h_0(n) = \LL_{\rm ref} h(n) = -g(n) \le -\delta h_0(n),
$$
where \beq\label{thm33f} \delta:=\inf_{n\in\nn^*}
\frac{g(n)}{h_0(n)} =\frac{2J_*}{h_0(N)}=\frac{2J_*}{\sum_{j=1}^N j
\frac{ \nu[j,N]}{\nu(j)}}. \nneq By Theorem \ref{thm-comparison},
for the metric $d(x,y):=h_0\circ d_G(x,y)$, $Ric(\LL, d)\ge \delta$.
Since
$$
n\le h_0(n)\le h_0(1) n =\frac{\nu[1,N]}{\nu(1)} n,
$$
by Theorem \ref{thm21}, we obtain
$$\aligned
W_{1,d_G}(P_t(x,\cdot), P_t(y,\cdot))&\le W_{1,d}(P_t(x,\cdot),
P_t(y,\cdot))\\
&\le e^{-\delta t} h_0(d_G(x,y))\\
&\le \frac{\nu[1,N]}{\nu(1)} e^{-\delta t} d_G(x,y),
\endaligned $$
i.e. the desired exponential convergence holds with $\delta$ given
in (\ref{thm33f}) and $K=\frac{\nu[1,N]}{\nu(1)} $, where $\nu$ is
given in (\ref{thm33e}). \nprf

\subsection{A discrete stochastic difference equation
in $\zz^d$}

\bexa\label{exa-sde} {\bf (A discrete stochastic difference equation
in $\zz^d$)} {\rm Let $S=\zz^d$ ($d\ge1$) equipped with the graph
metric $d_G(x,y)=\sum_{i=1}^d|x_i-y_i|$ (the $L^1$-metric). Given a
discrete vector field $b(x)=(b_1(x),\cdots,b_d(x)): \zz^d\to \rr^d$,
consider the generator \beq\label{prop3-10a} \LL f(x)=
\sum_{i=1}^d\left\{ \frac a2 [f(x+e_i)+f(x-e_i)-2f(x)]  + |b_i(x)|
[f(x+\sgn(b_i(x)) e_i )-f(x)]\right\}, \nneq where $a>0$ is a positive
constant, and $\sgn(r)=1, 0, -1$ according to $r>0, r=0, r<0$ (sign
of a real number $r$), and $(e_i)_j=\delta_{ij}$. Notice that
$$\sum_{i=1}^d [f(x+e_i)+f(x-e_i)-2f(x)]$$ is  the usual
discrete Laplacian operator on $\zz^d$ which is $2d$ times the graph
Laplacian $\Delta$ on $\zz^d$ used in this paper.

\bprop\label{prop3-10} If the discrete vector field $b(x)$ satisfies
the following dissipative condition at infinity w.r.t. the graph
metric $d_G$: \beq\label{prop3-10b} \sum_{i=1}^d
[(b_i(y)-b_i(x))\sgn(y_i-x_i) + |b_i(y)-b_i(x)| 1_{x_i=y_i}]\le -
\rho(d_{G}(x,y)), \nneq
where $\rho:\nn\to \rr$ is non-decreasing and
$$
\rho(k)=c(k-N), \ \forall k\ge N
$$
for some constants $c>0, N\in \nn$, then there are two constants
$K\ge 1, \kappa>0$ explicitly computable such that
$$
W_{1,d_G}(P_t(x,\cdot), P_t(y,\cdot))\le K e^{-\kappa t} d_G(x,y),\
t\ge0, x,y\in\zz^d.
$$
\nprop

The condition (\ref{prop3-10b}) can be viewed as the dissipativity at infinity w.r.t. the $L^1$-metric $d_{L^1}(x,y)=\sum_i |x_i-y_i|$, a counterpart of (\ref{dissipativity3}) w.r.t. the Euclidean metric, and
then this result is a lattice valued version of the $W_1$-exponential convergence obtained in Eberle \cite{Eberle16}, Luo and Wang \cite{LuoWang16} for diffusions.

\bprf {\bf Step 1. Construction of a coupling generator.} Given
$x\ne y$ in $\zz^d$ and $i=1,\cdots, d$, we want to construct a
coupling operator of $\LL_{i,x}$ and $\LL_{i,y}$ (acting only on the
$i$-th coordinate,  with $(x_j, y_j)_{j\ne i}$ fixed), where
$$
\LL_{i, x}f(x)=\frac a2 [f(x+e_i)+f(x-e_i)-2f(x)] + |b_i(x)|
[f(x+\sgn(b_i(x)) e_i)-f(x)].
$$

If $x_i\ne y_i$, let $\LL_i^\pi$ be the independent coupling of
$\LL_{i,x}$ and $\LL_{i,y}$.

If $x_i=y_i$,
$$\aligned
\LL_i^\pi F(x,y) &:= \frac a2[F(x+e_i,y+e_i)+ F(x-e_i, y-e_i)-2
F(x,y)] \\
&+ 1_{b_i(x)b_i(y)\ge 0} (|b_i(x)|\wedge
|b_i(y)|)[F(x+\sgn(b_i(x))e_i,y+ \sgn(b_i(y)) e_i)-F(x,y)]\\
&+1_{b_i(x)b_i(y)\ge 0} [|b_i(x)|- |b_i(x)|\wedge |b_i(y)| ]\cdot
[F(x+ \sgn(b_i(x))e_i, y) -F(x,y)]\\
&+1_{b_i(x)b_i(y)\ge 0} [|b_i(y)|- |b_i(x)|\wedge |b_i(y)| ]\cdot
[F(x, y+ \sgn(b_i(y))e_i) -F(x,y)]\\
&+1_{b_i(x)b_i(y)<0} |b_i(x)| \cdot
[F(x+\sgn(b_i(x)) e_i,y) -F(x,y)]\\
&+1_{b_i(x)b_i(y)<0} |b_i(y)| \cdot
[F(x, y+\sgn(b_i(y)) e_i) -F(x,y)].\\
\endaligned$$
The first term above means that the evolution of $(X_{t,i},
Y_{t,i})$ related to $a \Delta_i$ goes together once if $x_i=y_i$, as
well as the second term above related to the component $b_i$ and
that is with the maximal possible rate $|b_i(x)|\wedge |b_i(y)|$
once if $b_i(x)$ and $b_i(y)$ are of the same sign. The remained four
terms are related to free jumps.

Now define our  coupling generator of $\LL$ as
$$
\LL^\pi F(x,y)=\sum_{i=1}^d \LL_i^\pi F(x,y)
$$
and let $J^\pi((x,y),(x',y'))$ be the corresponding kernel of
jumping rates. By the construction above, if $d_G(x,y)=n\ge1$,
\beq\label{prop3-10c}\aligned
 a_\pi(x,y)&:=\sum_{d_G(x',y')=n-1}
J^\pi((x,y),(x',y'))\\
&= \sum_{i=1}^d 1_{x_i\ne y_i} \left\{a +
1_{b_i(x)\sgn(y_i-x_i)>0}|b_i(x)| +
1_{b_i(y)\sgn(y_i-x_i)<0}|b_i(y)| \right\};\\
b_\pi(x,y)&:=\sum_{d_G(x',y')=n+1}
J^\pi((x,y),(x',y'))\\
&=\sum_{i=1}^d 1_{x_i\ne y_i} \left\{a +
1_{b_i(x)\sgn(y_i-x_i)<0}|b_i(x)| +
1_{b_i(y)\sgn(y_i-x_i)>0}|b_i(y)| \right\}\\
&\quad+\sum_{i=1}^d 1_{x_i=y_i}\left\{ 1_{b_i(x)b_j(x)\ge0}
(|b_i(x)|+|b_i(y)|- 2 [|b_i(x)|\wedge
|b_i(y)|]\right.\\
&\quad\ \quad \quad \left. +  1_{b_i(x)b_j(x)<0} (|b_i(x)|+|b_i(y)|) ) \right\}\\
&=\sum_{i=1}^d 1_{x_i\ne y_i} \left\{a +
1_{b_i(x)\sgn(y_i-x_i)<0}|b_i(x)| +
1_{b_i(y)\sgn(y_i-x_i)>0}|b_i(y)| \right\}\\
&\quad+\sum_{i=1}^d 1_{x_i=y_i} |b_i(x)-b_i(y)|.
\endaligned \nneq
Therefore $\LL^\pi d_G(x,y)=b_\pi(x,y)-a_\pi(x,y)$ which is equal to
the left hand side of (\ref{prop3-10b}), and then $\le
-\rho(d_G(x,y))$.

Set $\beta_{-1}(x,y)=a_\pi(x,y)-a-\rho^+(d_G(x,y))$. It is
nonnegative: this is evident if $\rho(d_G(x,y))<0$, and if
$\rho(d_G(x,y))\ge0$,
$$
a_\pi(x,y)-a\ge b_\pi(x,y)-a + \rho(d_G(x,y))\ge \rho(d_G(x,y)).
$$
Then $a_\pi(x,y)=\beta_{-1}(x,y) + a + \rho^+(d_G(x,y))$. On the
other hand,
$$\aligned
b_\pi(x,y)&=\beta_{-1}(x,y) +a+\rho^+(d_G(x,y)) +(b_\pi(x,y)-a_\pi(x,y))\\
& \le \beta_{-1}(x,y) +a+\rho^+(d_G(x,y)) -\rho(d_G(x,y))\\
&= \beta_{-1}(x,y) +a+\rho^-(d_G(x,y)).
\endaligned
$$
In summary the comparison condition ${\bf C(J,\alpha,\beta)}$ is
satisfied for $\alpha=1$, $\beta_{1}=\beta_{-1}$ given above and
\beq\label{sde4} J_n(n-1)=a+\rho^+(n),\ J_n(n+1)=a+\rho^-(n), \
n\in\nn^*. \nneq

{\bf Step 2.} Let $h(n)=n$ for $n\ge N$. We have for any $n\ge N+1$,
$$\aligned
\LL_{\rm ref} h(n) &= J_n(n+1) D_+h(n+1) - J_n(n-1)
D_+h(n-1)\\
&=J_n(n+1)- J_n(n-1)=-\rho(n)=-c (n-N). \endaligned$$ Consider the
function $g(k)= c 1_{[0,N]}(k) + c 1_{k\ge N+1} (k-N)$. We will
construct $h(k)$ for $k\in [0, N-1]\cap \nn$ by solving
\beq\label{sde5} -\LL_{\rm ref} h(k)= J_k(k-1) D_+h(k-1) - J_k(k+1)
D_+h(k) =c, \ k\in [1,N]. \nneq Once it is found, we would get
$\LL_{\rm ref} h=-g$ on $\nn^*$.

Let $\nu$ be the symmetric measure of $\LL_{\mbox{ref}}$, determined by
$\nu(N)=1$ and $\nu(k)J_k(k+1)=\nu(k+1)J_{k+1}(k),  k\in\nn^*$ (the
detailed balance condition). Then for $n\in [1,N-1]$, \beq\label{nu}
\nu(n)=\prod_{k=n}^{N-1}
\frac{J_{k+1}(k)}{J_k(k+1)}=\frac{a^{N-n}}{\prod_{k=n}^{N-1}(a+\rho^-(k))
}. \nneq

 Let
$w(0,1)=\nu(1)J_1(0)$ and $w(k,k+1)= \nu(k)J_k(k+1)$ ($k\ge1$) (a
weight assigned to the edge $(k,k+1)$). Multiplying both sides of
(\ref{sde5}) by $\nu(k)$, we get
$$
w(k-1,k)D_+h(k-1) - w(k,k+1)D_+h(k)=c\nu(k),\ k\in [1,N].
$$
Summing this equality from $k=n$ to $N-1$, we get for all $1\le n\le
N-1$,
$$
w(n-1,n) D_+h(n-1) - w(N-1,N) D_+h(N-1) = c\nu[n, N-1].
$$
It remains to determine $D_+h(N-1)$. By (\ref{sde5}) for $k=N$, we
get by recalling that $h(n)=n$ for $n\ge N$ and (\ref{sde4}),
$$\aligned
w(N-1,N) D_+h(N-1)&=J_N(N-1)D_+h(N-1)\\
&=J_N(N+1) D_+h(N) + c=a+c. \endaligned$$ Therefore \beq\label{sde6}
D_+h(n-1)=\frac{c\nu[n,N] + a}{w(n-1,n)}, \ 1\le n\le N. \nneq As
for $n\le N$,
$$
\frac{w(n, n+1)}{w(n-1,n)}=
\frac{J_n(n+1)}{J_n(n-1)}=\frac{a+\rho^-(n)}{a}\ge 1,
$$
$D_+h(n-1)$ is decreasing in $n\le N$.
$$
D_+h(N-1)=
\frac{a+c}{w(N-1,N)}=\frac{a+c}{J_N(N-1)}=\frac{a+c}{a}\ge 1=
D_+h(N),
$$
$D_+h(n-1)$ is non-increasing over $\nn^*$.

Moreover
$$
h(n)-h(0)=\sum_{k=1}^n D_+h(k-1)= \sum_{k=1}^n \frac{c\nu[k,N] +
a}{w(k-1,k)}, \ n\in [1,N].
$$
Setting \beq\label{sde7}\aligned h_0(n)&:= h(n)-h(0)\\
\kappa&:=\inf_{k\ge 1} \frac{g(k)}{h_0(k)} = c\left(1+ \sum_{k=1}^N
\frac{c\nu[k,N] + a}{w(k-1,k)} \right)^{-1} \endaligned\nneq (the
last equality is obtained by calculus), we have for any $n\ge 1$,
$$
\LL_{\rm ref} h_0(n) =\LL_{\rm ref} h(n)=-g(n) \le -\kappa h_0(n).
$$
Applying Theorem \ref{thm-comparison}(a), $Ric(\LL, d)\ge\kappa$ for
$d(x,y)=h_0\circ d_G(x,y)$. As $D^+h$ is non-increasing
$$n=n
D_+h(n+N)\le h_0(n)=\sum_{k=1}^n D_+h(k)\le h_0(1) n$$ and
$h_0(1)=\frac{c\nu[1,N] + a}{\nu(1)a}$, we get by Theorem
\ref{thm21},
$$\aligned
W_{1,d_G}(P_t(x,\cdot), P_t(y,\cdot))&\le W_{1,d}(P_t(x,\cdot),
P_t(y,\cdot))\\
&\le e^{-\kappa t} h_0(d_G(x,y))\\
&\le \frac{c\nu[1,N] + a}{\nu(1)a} e^{-\kappa t} d_G(x,y).
\endaligned $$
\nprf } \nexa

\brmk {\rm
The coupling $\LL^\pi$ constructed in the proof of Proposition \ref{prop3-10} is optimal, i.e.
\beq\label{equiv}
Ric_{(x,y)} (\LL,d_G) =-\frac{\LL^\pi d_G(x,y)}{ d_G(x,y)}.
\nneq

In fact by the definition \ref{def_Ricci}, we have
\beq\label{left}
Ric_{(x,y)} (\LL,d_G) \ge -\frac{\LL^\pi d_G(x,y)}{ d_G(x,y)}.
\nneq
Now we turn to the "$\le$ " part. For every fixed $ x^0, y^0 \in S$, consider the following function $h_{(x^0, y^0)}: S \to \rr$,
$$
\begin{aligned}
h_{(x^0, y^0)} (z) &=\sum_{i=1}^d 1_{\{x_i^0 < y_i^0\}} z_i +\sum_{i=1}^d 1_{\{x_i^0 > y_i^0\}}( -z_i )\\
&\quad + \sum_{i=1}^d 1_{\{x_i^0 = y_i^0\}} \big( 1_{\{b_i(x^0) < b_i(y^0)\}} z_i + 1_{\{b_i(x^0) > b_i(y^0)\}}( -z_i ) \big) .
\end{aligned}
$$
It is easy to see that $$ \|h_{(x^0, y^0)}\|_{Lip(d_G)}= \sup_{x \neq y} \frac {|h_{(x^0, y^0)}(y)-h_{(x^0, y^0)}(x)|}{ d_G(x,y)} \le 1,$$ then by Kantorovich duality and definition \ref{def_Ricci}, for every fixed $x,y \in S$, we have by setting $J(x,x)=\lambda(y)$, $J(y,y)=\lambda(x)$,
\beq\label{right}
\begin{aligned}
&\quad Ric_{(x,y)} (\LL,d_G) \cdot d_G(x,y) \\
&= (\lambda (x)+ \lambda(y)) d_G(x,y)- T_{d_G(x,y)} (J(x,\cdot), J(y,\cdot))\\
& \le (\lambda (x)+ \lambda(y)) d_G(x,y)-\left( \sum_{z\in S} h_{(x, y)}(z)J(y,z)- \sum_{z\in S} h_{(x, y)}(z)J(x,z)\right)\\
&=-\sum_{i=1}^d[(b_i(y)-b_i(x))\sgn(y_i-x_i) + |b_i(y)-b_i(x)| 1_{x_i=y_i}]\\
&=-\LL^\pi d_G(x,y).
\end{aligned}
\nneq
\eqref{left} together with \eqref{right} implies \eqref{equiv}. }
\nrmk

\section{Zhong-Yang's estimate on graphs of nonnegative curvature}

On a Riemannian manifold $M$ of dimension $n$ without boundary or
with convex boundary $\partial D$, of bounded diameter $D\in
(0,+\infty)$ such that $Ric_x\ge 0$, the famous Zhong-Yang's
estimate \cite{ZY84} for the spectral gap $\lambda_1(\Delta_M)$ of the Laplace
operator $\Delta_M$ with the Neumann-boundary condition at $\partial
M$ (if it is not empty) says that
 \beq\label{zhongeyanges}
\lambda_1(\Delta_M)\ge \frac{\pi^2}{D^2}.
\nneq
The quantity $\frac{\pi^2}{D^2}$ is $\lambda_1([0,D])$ (the spectral gap of $\Delta$ with the Neumann boundary condition at the boundary $\{0,D\}$),
and $\eqref{zhongeyanges}$ becomes equality for
the circle $M=S^1$.

The following is a partial counterpart of Zhong-Yang's estimate on
graph.

\bthm\label{thm41} Assume that $Ric(\LL, d_G)\ge0$ and the
diameter $D_G$ of $(S,E)$ is finite. Let $\LL^\pi$ be a
$d_G$-optimal coupling generator, i.e.
$$
\LL^\pi d_G(x,y) = -Ric_{(x,y)}(\LL, d_G)\cdot d_G(x,y), (x,y) \in S^2.
$$
If for some positive constant $a>0$,
 \beq\label{thm41a} a_{\pi,1}(x,y) + 2
a_{\pi,2}(x,y)\ge
\begin{cases} a,\ &\text{ if } 1\le d_G(x,y)\le D_G-1;\\
2a, &\text{ if } d_G(x,y)=D_G,
\end{cases}
\nneq then for the metric
$$d(x,y):=\sin\frac{d_G(x,y)\pi}{2D_G},$$

 \beq\label{thm41b}
Ric(\LL, d)\ge 2a \left(1-\cos \frac{\pi}{2D_G}\right). \nneq
In particular for any eigenvalue $\lambda\ne 0$ of $-\LL$, $\displaystyle Re (\lambda)\ge Ric(\LL, d)\ge 2a \left(1-\cos \frac{\pi}{2D_G}\right). $
\nthm

\bprf By Lemma \ref{lem33} and its proof, we can assume that our $d_G$-optimal
coupling generator $\LL^\pi$ satisfies $
a_{\pi,2}(x,y)=b_{\pi,2}(x,y)=0 $. In that case the condition
(\ref{thm41a}) becomes
$$
a_{\pi,1}(x,y) \ge  a\ \text{ if }\ d_G(x,y)\le D_G-1;\ a_{\pi,1}(x,y) \ge 2a \text{ if
}\ d_G(x,y)=D_G.
$$
Let
$$
J_n(n-1)=J_n(n+1)=a, \ n\in [1,D_G-1], \ J_n(n-1)=2a, \ \text{ if }\
n=D_G,
$$
and $\beta_{-1}(x,y):= a_{\pi,1}(x,y)-J_n(n-1)$, where
$n=d_G(x,y)\ge1$. Since $Ric(\LL,d_G)\ge 0$, we have
$$
b_{\pi,1}(x,y)\le a_{\pi,1}(x,y)=\beta_{-1}(x,y)+J_n(n-1),$$ i.e. the
comparison condition $C(J;\alpha,\beta)$ is satisfied for $\alpha=1$
and $J,\beta$ given above. By direct calculus, for $h(k)=\sin
\frac{k\pi}{2D_G}$,
$$
\LL_{\rm ref} h(k) =- 2a \left(1-\cos \frac{\pi}{2D_G}\right) h(k), \
1\le k\le D_G.
$$
Therefore we get the desired result (\ref{thm41b}) by Theorem
\ref{thm-comparison}(a).
 \nprf

\brmk{\rm  Let us explain what our extra-condition \eqref{thm41a} means. At first a such type condition is indispensable, for $Ric(\vep \LL, d_G)\ge 0$ for any $\vep>0$, and $Ric(\vep \LL, d)=\vep Ric(\LL,d)$: so one requires some condition to specify this activity parameter $\vep$. Secondly our condition  \eqref{thm41a} is on the coupling generator: that is quite natural because the Ricci curvature is defined in that way.

Notice that if a coupling generator $\LL^\pi$ satisfies $\LL^\pi d(x,y)\le 0$ and
$$
[a_{\pi,1}(x,y)+ b_{\pi,1}(x,y) ] + 2 [a_{\pi,2}(x,y)+ b_{\pi,2}(x,y) ]\ge 2a
$$
(this can be interpreted as a $L^1$-volatility of $d(X_t,Y_t)$), then  \eqref{thm41a} is verified.
}\nrmk

\brmk{\rm  Since $1-\cos \frac{\pi}{2D_G}=\lambda_1([0,2D_G]\cap\nn)$, the spectral gap of the random walk on  $[0,2D_G]\cap\nn$ as seen from Example \ref{random walk},
besides the activity constant $a$ specified by \eqref{thm41a}, this theorem can be regarded as a counterpart of Zhong-Yang's estimate.

An exact counterpart of Zhong-Yang's estimate, under the condition     \eqref{thm41a},  could be formulated as: in the symmetric case, $\lambda_1(\LL)\ge a \lambda_1(\zz/(n\zz))$ where $n=2D_G$ or $n=2D_G +1$.

Our sentiment is that increasing $a_{\pi,2}(x,y)$ and  $b_{\pi,2}(x,y)$ in the coupling generator will yield better estimate of the spectral gap $\lambda_1$. On the Riemannian manifolds, that is
the reflection coupling. But as on general graphes the reflection coupling does not exist, we are content of the use of the one-step coupling in Lemma \ref{lem33}.
}\nrmk

\brmk{\rm  For general comparison theorems on the spectral gap of an elliptic diffusion generator on Riemannian manifolds, the reader is referred to
Chen M.F. and Wang F.Y. \cite{CW94, CW97} under some mixing condition on the curvature and dissipativity, and to Bakry-Qian \cite{BQ00} under the curvature lower bound condition.

}\nrmk

\section{Lyapunov function method}
Generally speaking the method of Lyapunov function yields
qualitative results about positive or exponential recurrence etc.
(see Meyn and Treedie \cite{MT}). However  two explicit quantitative
estimates of the exponential convergence rate based on Lyapunov
functions are known. The first one is due to D. Bakry {et al.} \cite{BBCG08} in the symmetric case: they gave an explicit estimate of
the spectral gap using the Lyapunov function together with the
spectral gap of the reflected process in a bounded domain. The other
is due to M. Hairer and Mattingly  \cite{HM11}, who gave a quantitative version of
Harris' theorem about the exponential convergence of $P^n$ to the
invariant probability measure $\mu$, for
 a single transition probability kernel $P(x,y)$.
Their condition is a combination of Lyapunov function (for controlling
the rate of returning time to small set) and a minorization
condition on the small set.

But the minorization condition in \cite{HM11} becomes,
in the continuous time case, a hypothesis on the unknown semigroup
$P_t$. The objective of this section is to replace the minorization
condition by some suitable one based on the generator $\LL$.

\bthm\label{thm-Lyapunov} Assume that

\benu
\item
there is some (Lyapunov) function $V:S\to [1,+\infty)$ such that for
some positive constants $r,b$ and some finite subset $K$  of $S$
$$
\LL V(x) \le -r V(x) + b 1_K(x)
$$
and
$$
\underline{V}(K^c):=\inf_{x\in K^c} V(x)>\frac{b}{r}
$$
(this last condition holds automatically by choosing $K$ large
enough, if $V(x)$ tends to infinity as $d_G(x,o)\to\infty$, where $o$ is some fixed point in $S$);

\item there are some pseudo-metric $d_\pi$ on $S$ (i.e. satisfying all
axioms of a metric except $d_\pi(x,y)$ may be zero for two different
points) and a coupling Markov generator $\LL^\pi$ of $\LL$ such that
for some constant $C>0$,
$$d_\pi(x,y)\le C (V(x)+V(y))$$  and

\beq\label{finiteC1} \LL^\pi d_\pi(x,y)\le - 1_{K^2}(x,y), \ \forall
(x,y)\in S^2\backslash \triangle. \nneq \nenu

 Then for the cost-function
$$
d_\beta(x,y) = d_\pi(x,y) + \beta 1_{x\ne y}[V(x)+V(y)]
$$
with the parameter $\beta$ verifying $0 <\beta<\frac{1}{2b}$, there is
some positive constant $\kappa$ (explicitly computable) such that
$Ric(\LL, d_\beta)\ge \kappa$. \nthm

\bprf  Note that
$$\LL^\pi [1_{\vartriangle^c}(V\oplus V)](x,y)\le
\LL^{\pi}(V\oplus V)(x,y)=\LL V(x)+\LL V(y),\ \forall (x,y)\in
S^2\backslash \triangle.$$

At first for $(x,y)\in (K^c)^2$ with $x\ne y$, by the condition
(\ref{finiteC1}) on $d_\pi$, we have
$$\aligned
\LL^\pi d_\beta(x,y)&\le \beta(\LL V(x) +\LL V(y))\\
 & \le -\beta
r(V(x)+V(y))\\
&\le -d_\beta(x,y)\cdot r\inf_{(x,y)\in (K^c)^2\backslash \triangle}\frac{\beta(V(x)+V(y))}{d_\pi(x,y)+ \beta(V(x)+V(y))}.\\
\endaligned$$
Now if $x\in K$ and $y\in K^c$,
$$\aligned
\LL^\pi d_\beta(x,y)&= \LL^\pi d_\pi(x,y) +\beta\LL^\pi(1_{\triangle^c}(V\oplus V))(x,y)\\
 &\le \beta(\LL V(x) +\LL V(y))\\
 & \le \beta b -r \beta V(y)\\
&\le - d_\beta(x,y) \inf_{x\in K,y\in K^c} \left[\frac{r\beta V(y) -b\beta
}{d_\pi(x,y)+\beta(V(x)+V(y))}\right].
\endaligned$$

Finally for $(x,y)\in K^2$ such that $x\ne y$,  we have
$$\aligned
\LL^\pi d_\beta(x,y)&\le \LL^\pi d_\pi(x,y) + \beta\LL^\pi
[1_{\vartriangle^c}(V\oplus V)](x,y)\\
&\le -1 + 2\beta b\\
&\le - d_\beta(x,y)\inf_{(x,y)\in K^2\backslash
\triangle}\left[\frac{1-2\beta b}{d_\pi(x,y)+ \beta(V(x)+V(y))}
\right]
\endaligned$$
provided that $2b\beta<1$. Summarizing the results in the three
cases above, we see that $Ric(\LL, d_\beta)\ge \kappa$ with
\beq\label{thm51d} \aligned \kappa:=\min &\left\{\frac{ \beta [r
\underline{V}(K^c) -b]
}{(C+\beta)(\overline{V}(K)+\underline{V}(K^c))};\  \frac{1-2\beta
b}{2(C+ \beta)\overline{V}(K)}\right\},
\endaligned\nneq
where $\overline{V}(K)=\sup_{x\in K}V(x)$. This Ricci curvature
lower bound is positive by the conditions in the theorem. \nprf

\brmk {\rm Given $\LL^\pi$, the smallest one of the metrics $d_\pi$
satisfying (\ref{finiteC1}) is
$$
d_\pi^0(x,y):=\ee_{(x,y)} \int_0^{\tau_c} 1_{K^2}(X_t,Y_t)dt
$$
where $(X_t,Y_t)$ is the coupling process generated by $\LL^\pi$,
starting from $(x,y)$, and $\tau_c:=\inf\{t\ge0; X_t=Y_t\}$ is the
coupling time. Obviously $\LL^\pi d_{\pi}^0(x,y)=-1_{K^2\backslash
\triangle}(x,y)$ and by the strong Markov property this metric
satisfies
$$
d^0_\pi(x,y)\le \max_{(x',y')\in K^2} d^0_\pi(x',y'),\ (x,y)\in S^2\backslash \triangle.
$$
Thus $d^0_\pi$ is bounded once if $\ee_{(x,y)}\tau_c<+\infty$. The
problem for our condition (\ref{finiteC1}) is to bound
$d^0_{\pi}(x,y)$ by $C(V(x)+V(y))$ for some explicit constant $C>0$.
 } \nrmk

We see that the proof above does not depend on the nearest-neighbor
condition for $J$. The following corollary is the counterpart of
Hairer-Mattingly's result in the continuous time case:

\bcor Without the nearest-neighbor condition for $J$ (i.e. $J(x,y)$
may be positive even if $d_G(x,y)\ge 2$), assume the Lyapunov
function condition in Theorem \ref{thm-Lyapunov}. If \beq \sum_{y\in
S} J(x_1,y)\wedge J(x_2,y)\ge \delta>0, \ x_1,x_2\in K, \nneq then
$d_\pi(x,y)=\frac 1{\delta}1_{x\ne y}$ satisfies (\ref{finiteC1}),
and the conclusion of Theorem \ref{thm-Lyapunov} holds. \ncor

\bprf For the discrete metric $d(x,y)=1_{x\ne y}$, by the proof of
Corollary \ref{cor28},
$$
Ric_{(x_1,x_2)} (\LL, d)= \sum_{y\in S} J(x_1,y)\wedge J(x_2,y),
$$
where $J(x_1,x_1)=\lambda(x_2), J(x_2,x_2)=\lambda(x_1)$, i.e. there is some coupling generator $\LL^\pi$ such that
$$
\LL^\pi d(x_1,x_2) \le -\delta, \ (x_1,x_2)\in K^2\backslash
\triangle.
$$
Then $d_\pi(x,y)=\frac 1{\delta}1_{x\ne y}$ satisfies
(\ref{finiteC1}) with $C=1/\delta$.
 \nprf

If $J$ is of nearest-neighbor type, to find $d_\pi$ satisfying
(\ref{finiteC1}), we can test
$$
d_\pi (x,y)= h\circ d_G(x,y), \ or\ d_\pi(x,y) = g(N-[d_G(x, K^c)
\wedge d_G(y, K^c)]),
$$
where $h,g:\nn\to\rr^+$ are nondecreasing functions such that
$h(n)=h(N),\ g(n)=g(N)$ for all $n\ge N$. That is the purpose of
the following corollary, which is a combination of the Lyapunov function and the Ricci curvature.

\bcor If $Ric(\LL,d_G)$ is bounded from below and
$J_*=\inf_{(x,y)\in E} J(x,y)>0$, then there exists a coupling
generator $\LL^\pi$ such that for each $N\ge1$, there is some
increasing function $h_0:\nn\to \rr$ with $h_0(0)=0$  such that
$D_+h_0$ is non-increasing and
$$
\LL^\pi h_0\circ d_G(x,y) \le - 1_{[1,N]} (d_G(x,y)).
$$
Therefore $d_\pi(x,y)=h_0\circ d_G(x,y)$ satisfies (\ref{finiteC1})
once if $N\ge Diam(K,d_G):=\max_{x,y\in K} d_G(x,y)$, and the
conclusion of Theorem \ref{thm-Lyapunov} holds true under the Lyapunov function condition (1) there.  \ncor

\bprf Let $R\ge0$ such that $Ric(\LL,d_G)\ge -R$. Following the
proof of Theorem \ref{thm33}, we can find a $d_G$-optimal coupling
generator $\LL^\pi$ satisfying the comparison condition ${\bf
C(J,\alpha,\beta)}$ for $\alpha=1, \beta_{-2}=\beta_{2}=0$ and
$$
J_n(n-1)=2J_*, \ J_n(n+1)=2J_* + R n, \ J_n(n\pm 2)=0.
$$
For any $N\ge 1$ fixed, if $h(n)=h(N)$ for all $n\ge N$, then
$$
\LL_{\rm ref} h(N) = J_N(N-1)(h(N-1)-h(N)).
$$
The solution $h_0$ satisfying $h_0(0)=0$ and $h_0(n)=h_0(N)$ for all
$n\ge N$ of the Poisson equation $\LL_{\rm ref} h(n)=-1,\ n\in
[1,N]$ is determined by
$$
D_+h_0(n-1)=\frac{\nu[n,N]}{2J_*\nu(n)}, n\in [1,N],
$$
where $\nu$ is the symmetric measure of $\LL_{\rm ref}$. By Theorem
\ref{thm-comparison},
$$
\LL^\pi h_0\circ d_G(x,y)\le (\LL_{\rm ref} h_0)(d_G(x,y)) \le
-1_{[1,N]} (d_G(x,y)).
$$
That completes the proof.
 \nprf

\section{Glauber dynamics and Gibbs sampler}
\subsection{Ricci curvature lower bound of Glauber dynamics }
In this subsection we consider the Ricci curvature lower bound of Glauber dynamics for approximating Gibbs measures in high dimension. Consider a fixed finite subset $V$ of $ \zz^d$, $\GG$ is a family of subsets of $V$ such that $ \bigcup_{\Lambda\in \GG} \Lambda= V $, i.e. a covering of $V$. For each $i\in V$, $(S_i, E_i)$ is a graph equipped with the metric $d_i$. For any $ \Lambda \in \GG$, denote by $S_\Lambda=\Pi_{i \in \Lambda} S_i$ the product graph, by $d_\Lambda (x_\Lambda, y_\Lambda) = \sum_{i \in \Lambda}
d_i(x_i,y_i)$ the $L^1$-metric on $S_\Lambda$.
Throughout this section, we consider the product graph $S_V= \Pi_{i \in V} S_i$ equipped with the $L^1$-metric $ d_{L^1}(x,y)=\sum_{i\in V} d_i(x_i,y_i)$. \\

Consider the generator:
\beq\label{jilin1}
\LL f (x) =\sum_{\Lambda \in \GG} \LL_\Lambda f(x),\ \LL_\Lambda f(x):=\sum_{x_\Lambda^\prime \in S_\Lambda} \big(
f(x^{x_\Lambda^\prime}) - f(x)\big) J_\Lambda(x,x_\Lambda^\prime),
\nneq
where $x^{x_\Lambda^\prime}$ is the new configuration such that
$$
(x^{x_\Lambda^\prime})_j=\begin{cases}
x_j,&\mbox{ if }j\notin \Lambda;\\
x_j^\prime, &\mbox{ if } j \in \Lambda,
\end{cases}
$$
and $J_\Lambda(x,x_\Lambda^\prime):=J(x,x^{x_\Lambda^\prime})$ is the jumps rate from $x$ to
$x^{x_\Lambda^\prime}$. Let $\lambda_\Lambda(x)=
{\displaystyle\sum_{x_\Lambda^\prime \neq x_\Lambda}} J_\Lambda (x,
x_{\Lambda}^\prime)$.

We assume

{\bf (H1)}: There exists a constant $\kappa_0>0$ such that for any $\Lambda \in \GG$ and the boundary condition
$x_{\Lambda^c}$ fixed,
$$
Ric(\LL_{\Lambda}, d_\Lambda)\ge \kappa_0 \ \text{ on }\ S_\Lambda.
$$

{\bf (H2)}: For all $j \in V$, $\Lambda \in \GG$ and for $x,y\in S_V$ such that $x_\cdot =
y_\cdot$ except the site $j$, there exists some $M =M(x,y)
\ge \max \{ \lambda_\Lambda(x), \lambda_\Lambda(y)\}$ and some
constant $C_{\Lambda j }>0$ such that
$$
W_{1,d_\Lambda} (\bar{J}_\Lambda (x,\cdot),\bar{J}_\Lambda
(y,\cdot)) \le C_{\Lambda j} d_j(x_j,y_j)
$$
where
$$
\begin{aligned}
&\bar{J}_\Lambda(x, \cdot)= \sum_{x_\Lambda^\prime \neq x_\Lambda} J_\Lambda (x,x_\Lambda^\prime)+(M-\lambda_\Lambda(x)) \delta_{x_\Lambda}(\cdot),\\
&\bar{J}_\Lambda(y, \cdot)= \sum_{y_\Lambda^\prime \neq y_\Lambda}
J_\Lambda (y,y_\Lambda^\prime)+(M-\lambda_\Lambda(y))
\delta_{y_\Lambda}(\cdot).
\end{aligned}
$$
The following result generalizes the Ligget's $M$-$\vep$ theorem in \cite[Theorem 3.1]{Wu06}.

\bthm\label{thm51} Assume {\bf (H1), (H2)}. For any $j \in \Lambda$,
denote by $ N (j) = |\{ \Lambda\in \GG: j \in \Lambda\}|$. If there
exists a constant $\kappa>0$ such that
$$
\kappa_0 N(j) -\sum_{\Lambda: j \notin \Lambda} C_{\Lambda j} \ge
\kappa
$$
for all $j \in V$, we have
$$
Ric(\LL, d_{L^1}) \ge \kappa.
$$
\nthm \bprf By Theorem \ref{thm21}, it is enough to construct a
coupling $ \LL^\pi $ of $\LL$ such that
$$
\LL^\pi d_{L^1}(x,y)\le -\kappa d_{L^1} (x,y).
$$
Following the proof of Corollary \ref{cor24} and the definition
of $d_{L^1}$, we have only to prove it for $x,y\in S_\Lambda$ such that $x_\cdot =y_\cdot$ except
some single site $j$.

So fix $x$ and $y$ such that $x_i=y_i$, $\forall i\neq j $ and
$x_j\neq y_j$. Let $J_\Lambda^\pi((x,y),\cdot)$ be the optimal
coupling of $\bar{J}_\Lambda(x, \cdot),\ \bar{J}_\Lambda(y, \cdot)$
for the Wasserstein metric $W_{1,d_\Lambda}$ and consider the
corresponding coupling generator $\LL_\Lambda^\pi$ of $\LL_\Lambda$
:
$$
\LL_\Lambda^\pi F(x,y)=\sum_{(x_{\Lambda}^\prime,
y_{\Lambda}^\prime)\in S_\Lambda \times
S_\Lambda}\left(F(x^{x_\Lambda^\prime},y^{y_\Lambda^\prime})-F(x,y)\right)J_\Lambda^\pi((x,y),(x_\Lambda^\prime,y_\Lambda^\prime)).
$$
We define coupling generator of $\LL$ by
$$
\LL^\pi F(x,y)=\sum_{\Lambda \in \GG}\LL_{\Lambda}^\pi F(x,y).
$$
Let us estimate $\LL_{\Lambda}^\pi d_{L^1}(x,y)$. Our discussion
will be separated into two cases: $j \in \Lambda$ and $j \notin
\Lambda$.

{\bf Case 1}: $j\in \Lambda$, by Theorem \ref{thm21} and {\bf (H1)},
we have \beq \label{51} \LL_\Lambda^\pi d_{L^1}(x,y)=\LL_\Lambda^\pi
d_\Lambda (x,y) \le -\kappa_0 d_\Lambda (x,y)=-\kappa_0
d_j(x_j,y_j). \nneq

{\bf Case 2}: $j \notin \Lambda$, by {\bf (H2)}, we have
\beq\label{52}
\begin{aligned}
\LL_\Lambda^\pi d_{L^1} (x,y) &=\sum_{(x_{\Lambda}^\prime, y_{\Lambda}^\prime)\in S_\Lambda \times S_\Lambda}\big(d_{L^1}(x^{x_\Lambda^\prime},y^{y_\Lambda^\prime})-d_{L^1}(x,y)\big) J_\Lambda^\pi((x,y), (x_\Lambda^\prime,y_\Lambda^\prime))\\
&=\sum_{(x_\Lambda^\prime,y_\Lambda^\prime)\neq (x_\Lambda,y_\Lambda)}\big(d_\Lambda(x_\Lambda^\prime,y_\Lambda^\prime)+d_j(x_j,y_j)-d_j(x_j,y_j)\big) J_\Lambda^\pi((x,y), (x_\Lambda^\prime,y_\Lambda^\prime))\\
&=\sum_{(x_\Lambda^\prime,y_\Lambda^\prime)\neq (x_\Lambda,y_\Lambda)}d_\Lambda(x_\Lambda^\prime,y_\Lambda^\prime) J_\Lambda^\pi((x,y), (x_\Lambda^\prime,y_\Lambda^\prime))\\
&=W_{1,d_\Lambda} (\bar{J}_\Lambda (x,\cdot),\bar{J}_\Lambda (y,\cdot))\\
&\le C_{\Lambda j} d_j(x_j,y_j).
\end{aligned}
\nneq Combine \eqref{51} and \eqref{52}, we get
$$
\begin{aligned}
\LL^\pi d_{L^1}(x,y) &=\sum_{\Lambda \in \GG} \LL_\Lambda^\pi d_{L^1}(x,y)\\
&=\sum_{\Lambda: j \in \Lambda} \LL_\Lambda^\pi d_{L^1}(x,y)+\sum_{\Lambda: j \notin \Lambda} \LL_\Lambda^\pi d_{L^1}(x,y)\\
&\le -\sum_{\Lambda: j \in \Lambda} \kappa_0 d_j(x_j,y_j)+\sum_{\Lambda: j \notin \Lambda} C_{\Lambda j} d_j(x_j,y_j)\\
&=-(\kappa_0 N(j) - \sum_{\Lambda: j \notin \Lambda}C_{\Lambda j}) d_j(x_j,y_j)\\
&\le -\kappa d_{L^1}(x,y).
\end{aligned}
$$
That completes the proof of this theorem. \nprf

Specially, when $V=\{1,2,\cdots,N\}$, $\GG=\{\{i\}, i\in V\}$,
$S=\Pi_{i=1}^N S_i$, we have the following result:

\bcor Assume {\bf (H1)} for $\LL_{\{i\}}=\LL_i$ and {\bf (H2)} for
$\{J_i(x,x_i^\prime); i \in V\}$ with $C_{\{i\}j}=C_{ij}$. If $
{\displaystyle\kappa_0-\sum_{i; i\ne j} C_{i j}} \ge \kappa$ for all $j \in
V$, then we have
$$
Ric(\LL,d_{L^1}) \ge \kappa.
$$
\ncor
\subsection{Gibbs sampler under the Dobrushin uniqueness condition:}
Let $\mu$ be a Gibbs measure on $S_V$. Consider the generator on
$S_V$ given by
$$
\LL f(x)=\sum_{i \in V} \LL_i f(x)= \sum_{ i \in V}(\mu_i(f|x)-f(x)),
$$
where $\mu_i(x_i| x)$ be the conditional distribution of $x_i$
knowing $(x_j)_{j\neq i}$ under $\mu$. We assume that $\mu_i(x_i|x)>0$ for all $x\in S_V$ and $x_i\in S_i$. Obviously $Ric_{(x_i,y_i)}(\LL_i, d_{i})=1$, i.e. $\kappa_0=1$ in {\bf (H1)}.
If $S_i$ is equipped with the complete
graph structure (i.e. $(x_i,y_i)\in E_i$ for every pair of two vertices $x_i,y_i\in S_i$), $\LL$ is of the nearest-neighbor type.

 This is the Gibbs sampler in
the continuous time case. Introduce the Dobrushin interdependence
coefficients:
$$
C_{ij}:=\sup_{x=y \mbox{ off }j}\frac{W_{1,d_i} (
\mu_i(\cdot|x),\mu_i(\cdot|y)  )}{ d_j(x_j,y_j)}, \ \forall i,j \in
V.
$$
(Obviously $C_{ii}=0$). With $J_i(x,y)=\mu_i(y_i|x)$, it coincides with $C_{{i}j}$ defined in {\bf (H2)}. We obtain thus
the following result whose equivalent version in $W_1$-exponential convergence was obtained by the third named author \cite{Wu06}.

\bcor Assume  {\bf the Dobrushin uniqueness
conditon} (\cite{Dobrushin68}), i.e.  there exists a constant $0< \kappa\le 1$ such that
$$
\sum_{i} C_{ij} \le 1-\kappa
$$
for every $j$,
then $Ric (\LL,d_{L^1}) \ge \kappa.$ \ncor

The same result in the discrete time case was proved by Ollivier \cite{Oll09}.

\subsection{Block Gibbs samplers under the Dobrushin-Shlosman analyticity conditon.}
When $V=[-N, N]^d$, $j\in \zz^d$,
\beq\label{jilin2}\GG= \{( [-l,l]^d+j) \cap V:
([-l,l]^d+j)\cap V \neq \emptyset )\}\nneq consider the generator on
$S_V$ given by
$$
\LL f(x)=\sum_{\Lambda \in \GG} \LL_\Lambda f(x)= \sum_{ \Lambda \in
\GG}(\mu_\Lambda(f|x)-f(x)),
$$
where $\mu_\Lambda(dx_\Lambda| x)$ be the conditional distribution
of $x_\Lambda$ knowing $x_{\Lambda^c}$. Assume the
Dobrushin-Shlosman analyticity condition (\cite{DoSh85}): there exist some positive
constants $C$ and $\delta$  such that if  $x_\cdot =y_\cdot $
except site $j$,
\beq\label{DoShA}
W_{1,d_\Lambda} (\mu_\Lambda(\cdot | x), \mu_\Lambda (\cdot|y)) \le
C e^{-\delta d_G(j,\Lambda)} d_j(x_j,y_j),
\nneq
where $d_G(j, \Lambda)= \inf \{ d_G(j,i): i \in \Lambda\}$ is the
distance of $j$ from $\Lambda$ and $d_G(j,i)=\sum_{k=1}^d|j_k-i_k|$ is the graph metric of $\zz^d$.

When the interaction is of bounded range $R$, i.e.
$W_{1,d_\Lambda} (\mu_\Lambda(\cdot | x), \mu_\Lambda (\cdot|y))=0$ for all configurations $x,y$ such that $x_\cdot=y_\cdot$ except site $j$ with
$d_G(j,\Lambda)>R$, the Dobrushin-Shlosman analyticity condition is equivalent to say that $x\to \mu_{\Lambda}(\cdot|x)$ is uniformly Lipschitzian from $(S_{\zz^d}, d_{L^1})$ to the space $\MM_1(S_\Lambda)$ of probability measures on $S_\Lambda$ equipped with the $L^1$-Wasserstein metric $W_{1,d_{L^1}}$.

\bcor Assume the Dobrushin-Shlosman analyticity condition (\ref{DoShA}) holds,
then   for any $l\ge0$ big enough,  there is some constant
$\kappa>0$ such that for the covering $\GG$ given by (\ref{jilin2}),
\beq\label{DoShb}
(2l+1)^d-\sum_{\Lambda\in \GG: j\notin \Lambda } C_{\Lambda j} \ge
\kappa.
\nneq
In that case we have
$$
Ric (\LL ,d_{L^1}) \ge \kappa.
$$
\ncor

\bprf In this case, $N(j) = (2l+1)^d$. The assumption {\bf (H1)} for
$\LL_\Lambda f(x)=\mu_\Lambda(f|x)-f(x)$ is satisfied for $\kappa_0
=1$. We note that {\bf (H2)} holds with $C_{\Lambda j} = C
e^{-\delta d(j,\Lambda)}$. Since the number of $\Lambda \in \GG$
such that $ d(j,\Lambda)=k$ is not larger than $2^d k^{d-1}
(2l+1)^{d-1}$, then there exists
some positive constant $ C^\prime$
$$
\begin{aligned}
\sum_{\Lambda: j \notin \Lambda } C_{\Lambda j} &\le C \sum_{\Lambda: j \notin \Lambda } e^{-\delta d(j,\Lambda)}\\
&\le C 2^d (2l+1)^{d-1}\sum_{k=1}^\infty e^{-\delta k} k^{d-1} \\
&= C^\prime  (2l+1)^{d-1}.
\end{aligned}
$$
Hence the condition (\ref{DoShb}) is satisfied with $\kappa=(2l+1-C')(2l+1)^{d-1}$ once  if $l>(C'-1)/2$. \nprf

The reader is referred to the lectures of Martinelli \cite{Martinelli} at Saint-Flour on Glauber dynamics (including the block type Gibbs samplers above) for the history and huge references on this topic.

\subsection{An interacting queue system}
\bexa[Queue system]{\rm Let $S = \nn^N$, equipped with the metric
$d_{L^1} (x,y) =\sum_{i=1}^N d_i(x_i,y_i)=\sum_{i=1}^N |x_i-y_i|$.
Consider the Gibbs measure
$$
\mu (x_1,x_2,\cdots,x_N) =\frac{e^{-\sum_{i\neq j} \beta_{i j} x_i
x_j } m(x_1) m(x_2) \cdots m(x_N)}{C},
$$
where $m(x_i)$ is the Poisson distribution of parameter $\lambda$,
$\beta_{i j}>0$ if $i \neq j$ ($\beta_{ii}=0$) are the correlation
coefficients of $x_i$ and $x_j$, $C$ is the normalization constant.
Then the conditional distribution of $x_i$ knowing $(x_j)_{j\neq i}$
is
$$
\mu_i(x_i| x)=\frac{e^{-\sum_{j:j\neq i} \beta_{ij}x_i x_j }
m(x_i)}{ \sum_{x_i} e^{-\sum_{j:j\neq i} \beta_{ij}x_i x_j }
m(x_i)}.
$$
Consider the Glauber dynamic:
$$
\begin{aligned}
\LL f(x) &= \sum_{i=1}^N \LL_i f(x)\\
&=\sum_{i=1}^N \left[  \lambda e^{-\sum_{j:j\neq i} \beta_{ij} x_j }
(f(x+e_i) - f(x)) + x_i ( f(x-e_i)- f(x) ) \right].
\end{aligned}
$$
If $1-\lambda \sum_{i:i\neq j} (1-e^{-\beta_{ij}})>0 $, we claim that
\beq\label{queue1} Ric(\LL, d_{L^1}) \ge 1-\lambda \sup_{j\in V}\sum_{i:i\neq j}
(1-e^{-\beta_{ij}}). \nneq

Indeed, let $J_i (x,x_i\pm 1):= J(x,x \pm e_i)$, in this model,
$$
\begin{aligned}
&J_i(x,x_i+1) =\lambda e^{-\sum_{j:j\neq i} \beta_{ij} x_j};\\
&J_i(x,x_i-1) = x_i.
\end{aligned}
$$
Consider $x_\cdot =y_\cdot $ except site $ j$ and $ y_j =x_j+1$.

For $i=j$, by Corollary \ref{cor2-10}, it is easy to see $ Ric
(\LL_i, d_i ) \ge 1$, then {\bf (H1)} holds with $\kappa_0 =1$.

For $i \neq j$, we use the following coupling:
$$
\begin{cases}
J_i^\pi ((x,y), (x_i+1,y_i+1) )=J_i (y,y_i+1)=\lambda e^{-(\sum_{k:k\neq i} \beta_{ik} x_k+ \beta_{ij})};\\
J_i^\pi ((x,y), (x_i-1,y_i-1) )=J_i (x,x_i-1)= x_i;\\
J_i^\pi ((x,y), (x_i+1,y_i) )=J_i (x,x_i+1)-J_i (y,y_i+1)=\lambda
e^{-\sum_{k:k\neq i} \beta_{ik} x_k}(1-e^{-\beta_{ij}}).
\end{cases}
$$
Then, we have
$$
\LL_i^\pi d_{L^1}(x,y) = \lambda e^{-\sum_{k:k\neq i} \beta_{ik}
x_k}(1-e^{-\beta_{ij}}) d_{L^1} (x,y).
$$
Thus  {\bf (H2)} holds with
$$
C_{ij}= \lambda (1-e^{-\beta_{ij}}).
$$
By Theorem \ref{thm51}, we get the result (\ref{queue1}).

The reader is referred to Dai Pra {\it et al.} \cite{DPP02} for the estimate of the spectral gap and of the rate in the exponential convergence in entropy for this model.
}
\nexa

\end{document}